\documentclass[a4paper]{article}
\usepackage{amsmath, amssymb, eucal, amscd, amstext}

\newtheorem{thm}{Theorem}[section]
\newtheorem{lem}[thm]{Lemma}
\newtheorem{eg}[thm]{Example}
\newtheorem{prop}[thm]{Proposition}
\newtheorem{cor}[thm]{Corollary}
\newtheorem{defn}[thm]{Definition}
\newtheorem{rem}[thm]{Remark}
\newtheorem{ntn}[thm]{Notation}

\newenvironment{prf}{{\noindent \textbf{Proof:}\ }}{\hfill $\Box$\\ \smallskip}

\numberwithin{equation}{section}

\newcommand{\abs}[1]{{\lvert#1\rvert}}

\newcommand{\id}{{\rm id}}
\newcommand{\ti}{\tilde}

\newcommand{\smnoind}{\smallskip\noindent}
\newcommand{\CL}{\mathcal{L}}

\newcommand{\sph}{\mathfrak{S}_1}

\newcommand{\BC}{\mathbb{C}}
\newcommand{\BN}{\mathbb{N}}
\newcommand{\BT}{\mathbb{T}}

\newcommand{\KF}{\mathfrak{F}}
\newcommand{\un}{{\rm unit}}
\newcommand{\BOI}{{\bigotimes}_{i\in I}}
\newcommand{\BOF}{{\bigotimes}_{i\in F}}
\newcommand{\tsi}{{\otimes}_{i\in I}\ \!}
\newcommand{\PI}{{\Pi}_{i\in I}}
\newcommand{\ut}{{\rm ut}}
\newcommand{\TBOIM}{{\ti\bigotimes}_{i\in I}}
\newcommand{\s}{{\rm s}}
\newcommand{\ct}{{\rm ct}}
\newcommand{\ii}{{i\in I}}

\begin{document}

\title{On genuine infinite algebraic tensor products \footnote{This work is supported by
the National Natural Science Foundation of China (10771126)}}
\author{Chi-Keung Ng}

\maketitle
\begin{abstract}
In this paper, we study genuine infinite tensor products of some algebraic structures. 
By a genuine infinite tensor product of vector spaces, we mean a vector space $\BOI X_i$ whose linear maps coincide with multilinear maps on an infinite family $\{X_i\}_\ii$ of vector spaces.
After establishing its existence, we give a direct sum decomposition of $\BOI X_i$ over a set $\Omega_{I;X}$, through which we obtain a more concrete description and some properties of $\BOI X_i$. 
If $\{A_i\}_{i\in I}$ is a family of unital $^*$-algebras, we define, through a subgroup $\Omega^\ut_{I;A}\subseteq \Omega_{I;A}$, an interesting subalgebra $\BOI^\ut A_i$. 
When all $A_i$ are $C^*$-algebras or group algebras, it is the linear span of the tensor products of unitary elements of $A_i$.
Moreover, it is shown that $\BOI^\ut \BC$ is the group algebra of $\Omega^\ut_{I;\BC}$. 
In general, $\BOI^\ut A_i$ can be identified with the algebraic crossed product of a cocycle twisted action of $\Omega^\ut_{I;A}$. 
On the other hand, if $\{H_i\}_{i\in I}$ is a family of inner-product spaces, we define a Hilbert $C^*(\Omega^\ut_{I;\BC})$-module $\bar\bigotimes^{\rm mod}_\ii H_i$, which is the completion of a subspace $\BOI^\un H_i$ of $\BOI H_i$. 
If $\chi_{\Omega^\ut_{I;\BC}}$ is the canonical tracial state on $C^*(\Omega^\ut_{I;\BC})$, then $\bar\bigotimes^{\rm mod}_\ii H_i\otimes_{\chi_{\Omega^\ut_{I;\BC}}}\BC$ coincides with the Hilbert space $\bar\bigotimes^{\phi_1}_\ii H_i$ given by a very elementary algebraical construction and is a natural dilation of the infinite direct product $\prod \tsi H_i$ as defined by J. von Neumann.  
We will show that the canonical representation of $\BOI^\ut \CL(H_i)$ on $\bar\bigotimes^{\phi_1}_\ii H_i$ is injective (note that the canonical representation of $\BOI^\ut \CL(H_i)$ on $\prod \tsi H_i$ is non-injective). 
We will also show that if $\{A_i\}_\ii$ is a family of unital Hilbert algebras, then so is $\BOI^\ut A_i$. 

\smnoind
MSC 2010: Primary: 15A69, 46M05; Secondary: 16G99, 16S35, 20C07, 46C05, 46L99, 47A80

\smnoind
Keywords: infinite tensor products; unital $^*$-algebras; twisted crossed products; inner product spaces; representations

\end{abstract}

\section{Introduction}

\medskip

In this paper, we study infinite tensor products of some algebraic structures. 
In the literature, infinite tensor products are often defined as inductive limit of finite tensor products (see e.g. \cite{Black}, \cite{Bru} \cite{Flor}, \cite{GN}, \cite{Pow}). 
As far as we know, the only alternative approach so far is the one by J. von Neumann, concerning \emph{infinite direct products of Hilbert spaces} (see \cite{vN}).
Some authors used this approach to define infinite tensor products of other functional analytic structures (see e.g. \cite{BC}, \cite{Gill} and \cite{Gui}). 
The work of von Neumann attracted the attention of many physicists who are interested in ``quantum mechanics with infinite degrees of freedom'', as well as mathematicians whose interest is in the field of operator algebras (see e.g. \cite{AN}, \cite{BeC}, \cite{BC}, \cite{EK}, \cite{GS}, \cite{Sto71}, \cite{TW}). 

\medskip

However, von Neumann's approach is not appropriate for purely algebraic objects. 
The aim of this article is to study ``genuine infinite algebraic tensor products'' (i.e. ones that are defined in terms of multilinear maps instead of through inductive limits) of some algebraic structures. 
There are several motivations behind this study.

\smnoind 1. Conceptually speaking, it is natural to define ``infinite tensor products'' as the object that produces a unique linear map from a multilinear map on a given infinite family of objects (see Definition \ref{def:inf-ten-prod}). 
As infinite direct products of Hilbert spaces are important in both Physics and Mathematics, it is believed that such infinite tensor products of algebraic structures are also important. 

\smnoind 2. We want to construct an infinite tensor product of Hilbert spaces that is easier for non-analyst to grasp (compare with the infinite direct product as defined by J. von Neumann; see Lemma \ref{lem:phi0-innpro} and Remark \ref{rem:cp-inf-direct-prod}(d)) and is more natural (see Theorem \ref{thm:ten-prod-hil-cT-mod}, Example \ref{eg:ten-Hil-mod} and Example \ref{eg:ten-prod-rep-BC}). 

\smnoind 3. Given a family of groups $\{G_i\}_{i\in I}$, it is well-known that the group algebra of the group
$${\bigoplus}_{i\in I} G_i := \big\{[g_i]_{i\in I}\in \PI G_i: g_i = e \text{ except for finite number of }i\in I\big\}$$ 
is an inductive limit of finite tensor products. 
However, if one wants to consider the group algebra $\BC[\PI G_i]$, one is forced to consider a ``bigger version of tensor products'' (see Example \ref{eg:gp-alg-inf-prod-gp}). 

\medskip

In this article, the algebraic structures that we concern with are vector spaces, unital $^*$-algebras, inner-product spaces as well as $^*$-representations of unital $^*$-algebras on Hilbert spaces. 
In our study, we discovered some interesting phenomena of infinite tensor products that do not have counterparts in the case of finite tensor products. 
Most of these phenomena related to certain object, $\Omega_{I;X}$, defined as in Remark \ref{rem:alt-con-ten-prod}(d), which ``encodes the asymptotic information'' of a given family $\{X_i\}_{i\in I}$. 

\medskip

In Section 2, we will begin our study by defining the infinite tensor product $(\BOI X_i, \Theta_X)$ of a family $\{X_i\}_{i\in I}$ of vector spaces. 
Two particular concerns are bases of $\BOI X_i$ as well as the relationship between $\BOI X_i$ and inductive limits of finite tensor products of $\{X_i\}_{i\in I}$ (which depend on choices of fixed elements in $\PI X_i$). 
In order to do these, we obtain a direct sum decomposition of $\BOI X_i$ indexed by a set $\Omega_{I;X}$ (see Theorem \ref{thm:bas-inf-ten}) with all the direct summand being inductive limits of finite tensor products (see Proposition \ref{prop:ju-inj}(b)). 
From this, we also know that the canonical map 
$$\Psi: \BOI L(X_i;Y_i) \to L(\BOI X_i;\BOI Y_i)$$ 
is injective (but not surjective). 
As a consequence, $\BOI X_i$ is automatically a faithful module over the big unital commutative algebra $\BOI \BC$ (see Corollary \ref{cor:ten-prod-as-mod} and Example \ref{eg:u-ten-prod}). 
Moreover, one may regard the canonical map 
$$\Theta_\BC : \PI \BC \to \BOI  \BC$$ 
as a generalised multiplication (see Example \ref{eg:u-ten-prod}(a)). 
In this sense, one can make sense of infinite products like $(-1)^I$. 

\medskip

Clearly, $\BOI A_i$ is a unital $^*$-algebra if all $A_i$ are unital $^*$-algebras. 
We will study in Section 3, a natural $^*$-subalgebra $\BOI^\ut A_i$ of $\BOI A_i$ which is a direct sum over a subgroup $\Omega^\ut_{I;A}$ of the semi-group $\Omega_{I;A}$. 
The reasons for considering this subalgebra are that it has good representations (see the discussion after Proposition \ref{prop:ten-prod-st-alg}), and it is 
big enough to contain $\BC[\PI G_i]$ when $A_i = \BC[G_i]$ for all $i\in I$ (see Example \ref{eg:gp-alg-inf-prod-gp}(a)). 
Moreover, if all $A_i$ are generated by their unitary elements (in particular, if $A_i$ are group algebras or $C^*$-algebras), then $\BOI^\ut A_i$ is the linear span of the tensor products of unitary elements in $A_i$.
We will show that $\BOI^\ut A_i$ can be identified with the crossed products of some twisted actions in the sense of Busby and Smith (i.e., a cocycle action with a $2$-cocycle) of $\Omega^\ut_{I;A}$ on $\BOI^e A_i$ (the unital $^*$-algebra inductive limit of finite tensor products of $A_i$). 
Moreover, it is shown that $\BOI^\ut \BC$ can be identified with the group algebra of $\Omega_{I;\BC}^\ut$ (Corollary \ref{cor:boiut-C}). 
We will also study the center of $\BOI^\ut A_i$ in the case when $A_i$ is generated by its unitary elements (for all $\ii$). 

\medskip

In Section 4, we will consider tensor products of inner-product spaces. 
If $\{H_i\}_{i\in I}$ is a family of inner-product spaces, we define a natural inner-product on a subspace $\BOI^\un H_i$ of $\BOI H_i$ (see Lemma \ref{lem:phi0-innpro}(b)). 
In the case of Hilbert spaces, the completion $\bar\bigotimes_{i\in I}^{\phi_1} H_i$ of $\BOI^\un H_i$ is a ``natural dilation'' of the infinite direct product $\prod \tsi H_i$ as defined by J. von Neumann in \cite{vN} (see Remark \ref{rem:cp-inf-direct-prod}(b)). 
Note that the construction for $\bar\bigotimes_{i\in I}^{\phi_1} H_i$ is totally algebraical and is more natural (see Example \ref{eg:ten-Hil-mod} and Example \ref{eg:ten-prod-rep-BC}). 
Note also that one can construct $\prod \tsi H_i$ in a similar way as $\bar\bigotimes_{i\in I}^{\phi_1} H_i$ (see Remark \ref{rem:cp-inf-direct-prod}(d)). 
On the other hand, there is an inner-product $\BC[\Omega^\ut_{I;\BC}]$-module structure on $\BOI^\un H_i$ which produces $\bar\bigotimes_{i\in I}^{\phi_1} H_i$ (see Theorem \ref{thm:ten-prod-hil-cT-mod}), as well as many other pre-inner-products on $\BOI^\un H_i$ (see Remark \ref{rem:hil-mod}(a)). 

\medskip

Section 5 will be devoted to the study of $^*$-representations of unital $^*$-algebras. 
More precisely, if $\Psi_i:A_i\to \CL(H_i)$ is a unital $^*$-representations ($i\in I$), we define a canonical $^*$-representation 
$$\BOI^{\phi_1} \Psi_i\ :\ \BOI^\ut A_i\ \to \ \CL\big(\bar\bigotimes_{i\in I}^{\phi_1} H_i\big).$$ 
We will show in Theorem \ref{thm:inf-ten-c-st-alg}(c) that if all $\Psi_i$ are injective, then $\BOI^{\phi_1} \Psi_i$ is also injective. 
This is equivalent to the canonical $^*$-representations of $\BOI^\ut \CL(H_i)$ on $\bar\bigotimes_{i\in I}^{\phi_1} H_i$ being injective, and is related to the ``strong faithfulness'' of the canonical action of $\Omega^\ut_{I;\CL(H)}$ on $\Omega^\un_{I;H}$ (see Remark \ref{rem-act-OA-OH}(b)). 
Note however, that the corresponding tensor type representation of $\BOI^\ut \CL(H_i)$ on $\prod \tsi H_i$ is non-injective. 
Consequently, if $(H_i,\pi_i)$ is a unitary representation of a group $G_i$ that induced an injective $^*$-representation of $\BC[G_i]$ on $H_i$ ($i\in I$), then we obtain injective ``tensor type'' $^*$-representation of $\BC[\PI G_i]$ on $\bar\bigotimes_{i\in I}^{\phi_1} H_i$ (see Corollary \ref{cor:ten-rep-prod-gps}).
On the other hard, we will show that $\bigoplus_{\rho\in \PI S(A_i)} \big( {\bar\bigotimes}_{i\in I}^{\phi_1} H_{\rho_i}, \BOI^{\phi_1} \pi_{\rho_i}\big)$ is an injective $^*$-representation of $\BOI^\ut A_i$ when all $A_i$ are $C^*$-algebras (Corollary \ref{cor:spat-ten-prod}). 
Finally, we show that if all $A_i$ are unital Hilbert algebras, then so is $\BOI^\ut A_i$. 

\medskip

\begin{ntn}
i). In this article, all the vector spaces, algebras as well as inner-product spaces are over the complex field $\BC$, although some results remain valid if one considers the real field instead.

\smnoind
ii). Throughout this article, $I$ is an infinite set, and $\KF$ is the set of all non-empty finite subsets of $I$. 

\smnoind
iii). For any vector space $X$, we write $X^\times := X\setminus \{0\}$ and put $X^*$ to be the set of linear functionals on $X$. 
If $Y$ is another vector space, we denote by $X\otimes Y$ and $L(X;Y)$ respectively, the algebraic tensor product of $X$ and $Y$, and the set of linear maps from $X$ to $Y$.
We also write $L(X) := L(X;X)$. 

\smnoind
iv). If $\{X_i\}_{i\in I}$ is a family of vector spaces and $x\in \PI X_i$, we denote by $x_i$ the ``$i^{\rm th}$-coordinate'' of $x$ (i.e. $x = [x_i]_{i\in I}$). 
If $x,y\in \PI X_i$ such that $x_i = y_i$ except for a finite number of $i\in I$, we write 
\begin{quotation}
$x_i = y_i\ $ e.f. 
\end{quotation}

\smnoind
v). If $V$ is a normed spaces, we denote by $\CL(V)$ and $V'$ the set of bounded linear operators and the set of bounded linear functionals respectively,  on $V$.
Moreover, we set $\sph(V):= \{x\in V: \|x\| =1\}$ as well as $B_1(V):= \{x\in V: \|x\| \leq 1\}$. 

\smnoind
vi). If $A$ is a unital $^*$-algebra, we denote by $e_A$ the identity of $A$ and $U_A:= \{a\in A: a^*a = e_A = aa^*\}$. 
\end{ntn}

\medskip

\section{Tensor products of vector spaces}

\medskip

\emph{In this section, $\{X_i\}_{i\in I}$ and $\{Y_i\}_{i\in I}$ are families of non-zero vector spaces.} 

\medskip

\begin{defn}\label{def:inf-ten-prod}
Let $Y$ be a vector space. A map $\Phi: \PI X_i\to Y$ is said to be \emph{multilinear} if $\Phi$ is linear on each variable. 
Suppose that $\bigotimes_{i\in I} X_i$ is a vector space and $\Theta_X: \PI X_i\to \bigotimes_{i\in I} X_i$ is a multilinear map such that for any vector space $Y$ and any multilinear map $\Phi: \PI X_i\to Y$, there exists a unique linear map $\tilde \Phi: \bigotimes_{i\in I} X_i \to Y$ with $\Phi = \tilde \Phi\circ \Theta_X$. 
Then $\left(\bigotimes_{i\in I} X_i, \Theta_X\right)$ is called the \emph{tensor product} of $\{X_i\}_{i\in I}$. 
We will denote $\tsi x_i := \Theta_X(x)$ ($x\in \PI  X_i$) and set $X^{\otimes I} := \bigotimes_{i\in I} X_i$ if all $X_i$ are equal to the same vector space $X$.
\end{defn}

\medskip

Let us first give the following simple example showing that non-trivial multilinear maps with infinite number of variables do exist. 
They are also crucial for some constructions later on. 

\medskip

\begin{eg}\label{eg:multi-linear}
(a) Let $\PI^1 \BC := \{\beta\in \PI \BC: \beta_i = 1 \text{ e.f.}\}$ and set 
$$\varphi_1(\beta)
\ :=\ \begin{cases}
\Pi_{i\in I} \beta_i \  &\text{if } \beta\in \PI^1 \BC\\
0 & \text{otherwise}.
\end{cases}$$ 
It is not hard to check that $\varphi_1$ is a non-zero multilinear map from $\PI \BC$ to $\BC$. 
If $\phi_1: \BOI\BC \to \BC$ is the linear functional induced by $\varphi_1$ (the existence of $\BOI \BC$ will be established in Proposition \ref{prop:exist-ten-prod}(a)), then $\phi_1$ is an involutive unital map. 

\smnoind
(b) Let $\PI^0 \BC := \{\beta\in \PI  \BC: \sum_{i\in I} \abs{\beta_i -1} < \infty\}$. 
For each $\beta\in \PI^0  \BC$, the net $\{\Pi_{i\in F} \beta_i\}_{F\in \KF}$ converges to a complex number, denoted by $\PI \beta_i$ (see e.g. \cite[2.4.1]{vN}). 
We define $\varphi_0(\beta) :=\ \Pi_{i\in I} \beta_i$ whenever $\beta\in \PI^0 \BC$ and set $\varphi_0|_{\PI\BC \setminus \PI^0 \BC} \equiv 0$. 
As in part (a), $\varphi_0$ induces an involutive unital linear functional $\phi_0$
on $\BOI\BC$. 
\end{eg}

\medskip

Clearly, infinite tensor products are unique (up to linear bijections) if they exist. 
The existence of infinite tensor products follows from a similar argument as that for finite tensor products, but we give an outline here for future reference. 

\medskip

\begin{prop}\label{prop:exist-ten-prod}
(a) The tensor product $\left(\bigotimes_{i\in I} X_i, \Theta_X\right)$ exists. 

\smallskip\noindent
(b) If $\{A_i\}_{i\in I}$ is a family of algebras (respectively, $^*$-algebras), then $\bigotimes_{i\in I} A_i$ is an algebra (respectively, a $^*$-algebra) with
$(\tsi a_i)(\tsi b_i) := \tsi a_ib_i$ (and $(\tsi a_i)^* := (\tsi a_i^*)$) for $a,b\in \PI  A_i$. 

\smnoind
(c) If $\Psi_i: A_i \to L(X_i)$ is a homomorphism for each $i\in I$, there is a canonical homomorphism $\TBOIM \Psi_i: \bigotimes_{i\in I} A_i \to L\left(\bigotimes_{i\in I} X_i\right)$ such that $\big(\TBOIM \Psi_i\big)(\tsi a_i)\tsi x_i = \tsi \Psi_i(a_i)x_i$ ($a\in \PI A_i$ and $x\in \PI X_i)$. 

\smnoind
(d) If $A = \bigoplus_{n=0}^\infty A_n$ is a graded algebra and $\bigoplus_{n=0}^\infty M_n$ is a graded left $A$-module, then $\bigoplus_{n=0}^\infty \bigotimes_{k\geq n} M_k$ is a graded $A$-module with
$a_m (\otimes_{k\geq n} x_k) = \otimes_{k\geq n} a_mx_k\in \bigotimes_{k \geq m+n} M_k$ ($a_m\in A_m; x\in \Pi_{k\geq n} M_k)$. 
\end{prop}
\begin{prf}
Parts (b), (c) and (d) follow from the universal property of tensor products, 
and we will only give a brief account for part (a).  
Let $V$ be the free vector space generated by elements in $\PI  X_i$ and $\Theta_0: \PI  X_i\to V$ be the canonical map. 
Suppose that $W:= {\rm span}\ \!W_0$, where
\begin{eqnarray}\label{eqt:def-W0}
\lefteqn{W_0 \ :=\ \big\{\lambda \Theta_0(u) + \Theta_0(v) - \Theta_0(w): \lambda\in \BC; u,v,w\in \PI  X_i; \exists i_0\in I \text{ with }} \nonumber\\ 
&& \qquad \qquad \qquad \qquad \qquad \qquad \qquad  \lambda u_{i_0} + v_{i_0} = w_{i_0} \text{ and } u_j = v_j = w_j, \forall j\in I\setminus\{i_0\}\big\}.
\end{eqnarray}
If we put $\bigotimes_{i\in I} X_i := V/W$, and set $\Theta_X$ to be the composition of $\Theta_0$ with the quotient map from $V$ to $V/W$, then they will satisfy the requirement in Definition \ref{def:inf-ten-prod}.
\end{prf}

\medskip

In the following remark, we list some observations that may be used implicitly throughout this article.

\medskip

\begin{rem}\label{rem:alt-con-ten-prod}
(a) As $\Theta_X$ is multilinear, $\BOI X_i = {\rm span}\ \! \Theta_X\big(\PI  X_i^\times \big)$. 

\smnoind
(b) If $I_1$ and $I_2$ are non-empty disjoint subsets of $I$ with $I = I_1\cup I_2$, it follows, from the universal property, that $\bigotimes_{i\in I} X_i \cong \big(\bigotimes_{i\in I_1} X_i \big) \otimes \big( \bigotimes_{j\in I_2} X_j \big)$ canonically. 

\smnoind
(c) $\BOI (X_i\otimes Y_i) \cong (\BOI X_i)\otimes (\BOI Y_i)$ canonically. 

\smnoind
(d) For any $x, y\in \PI  X_i^\times$, we denote 
\begin{equation*}
x \sim y \quad \text{if} \quad x_i = y_i \ \ e.f. 
\end{equation*}
Obviously, $\sim$ is an equivalence relation on $\PI  X_i^\times$, and we set $[x]_\sim$ to be the equivalence class of $x\in \PI  X_i^\times$. 
Let $\Omega_{I;X}$ be the collection of such equivalence classes. 
It is not hard to see that $\Omega_{I;\BC}$ is a quotient group of $\PI \BC^\times$, and that it acts freely on $\Omega_{I;X}$.  

\smnoind
(e) The element $\tsi 1\in \BC^{\otimes I}$ is non-zero. 
In fact, if $\tsi 1 = 0$, then $\BC^{\otimes I} = (0)$ (by Proposition \ref{prop:exist-ten-prod}(b)), and this implies the only multilinear map from $\PI \BC$ to $\BC$ being zero, which contradicts Example \ref{eg:multi-linear}. 
\end{rem}

\medskip

The ``asymptotic object'' $\Omega_{I;X}$ as defined in part (c) above is crucial in the study of genuine infinite tensor product, as can be seen in our next result. 
Let us first give some more notations here. 
For every $u\in \PI  X_i^\times$, we set 
$$\PI^u X_i := \{x\in \PI  X_i: x\sim u\}
\quad \text{and}\quad 
\BOI^u X_i := {\rm span}\ \! \Theta_X(\PI^u  X_i).$$ 
If $u\sim v$, then $\PI^u  X_i = \PI^v  X_i$, and we will also denote $\PI^{[u]_\sim}  X_i := \PI^{u}  X_i$ as well as $\bigotimes_{i\in I}^{[u]_\sim} X_i := \bigotimes_{i\in I}^{u} X_i$.

\medskip

\begin{thm}\label{thm:bas-inf-ten}
$\bigotimes_{i\in I} X_i = \bigoplus_{\omega\in \Omega_{I;X}} \bigotimes_{i\in I}^{\omega} X_i$.
\end{thm}
\begin{prf}
Suppose that $x^{(1)},...,x^{(n)}\in \PI  X_i^\times$ and $0 = n_0 < \cdots < n_N = n$ is a sequence satisfying
$x^{(n_k+1)} \sim \cdots \sim x^{(n_{k+1})}$ for $k \in  \{0,...,N-1\}$, but $x^{(n_k)}\nsim x^{(n_l)}$ whenever $1\leq k\neq l\leq N$. 
We first show that if $\nu_1,...,\nu_n\in \BC$ with $\sum_{l=1}^n \nu_l \Theta_X(x^{(l)}) = 0$, then 
$${\sum}_{l= n_k + 1}^{n_{k+1}} \nu_l \Theta_X(x^{(l)})\ =\ 0 \qquad  (k = 0,..., N-1).$$ 

In fact, by the proof of Proposition \ref{prop:exist-ten-prod}(a), there exist $m\in \BN$, $\mu_1,..., \mu_m\in \BC$ and $\lambda_k \Theta_0(u^{(k)}) + \Theta_0(v^{(k)}) - \Theta_0(w^{(k)})\in W_0$ ($k=1,...,m$) such that 
\begin{equation*}
{\sum}_{l=1}^n \nu_l\Theta_0(x^{(l)})
 = {\sum}_{k=1}^m \mu_k\big(\lambda_k \Theta_0(u^{(k)}) + \Theta_0(v^{(k)}) - \Theta_0(w^{(k)})\big).
\end{equation*}
Observe that if one of the elements in $\{ u^{(k)}, v^{(k)}, w^{(k)}\}$ is equivalent to $x^{(1)}$ (under $\sim$), then so are the other two (see \eqref{eqt:def-W0}). 
After renaming, one may assume that $u^{(k)}\sim v^{(k)}\sim w^{(k)}\sim x^{(1)}$ for $k=1,...,m_1$, but none of $u^{(k)}$, $v^{(k)}$ nor $w^{(k)}$ is equivalent to $x^{(1)}$ when $k\in \{m_1+1,...,m\}$. 

Since the two sets 
$$\{x^{(n_1+1)},...,x^{(n)}\}\cup \{u^{(m_1+1)},...,u^{(m)}\}\cup 
\{v^{(m_1+1)},...,v^{(m)}\}\cup \{w^{(m_1+1)},...,w^{(m)}\}$$
and $\{x^{(1)},...,x^{(n_1)}\}\cup \{u^{(1)},...,u^{(m_1)}\}\cup \{v^{(1)},...,v^{(m_1)}\}\cup \{w^{(1)},...,w^{(m_1)}\}$ are disjoint and elements in $\Theta_0\left( \PI  X_i\right)$ are linearly independent in $V$, we have 
\begin{equation*}
{\sum}_{l= 1}^{n_1} \nu_l \Theta_0(x^{(l)})
\ - \ {\sum}_{k=1}^{m_1} \mu_k\big(\lambda_k \Theta_0(u^{(k)}) + \Theta_0(v^{(k)}) - \Theta_0(w^{(k)})\big)
\ = \ 0.
\end{equation*}
This implies that $\sum_{l= 1}^{n_1} \nu_l \Theta_X(x^{(l)}) = 0$. 
Similarly, $\sum_{l= n_k+1}^{n_{k+1}} \nu_l \Theta_X(x^{(l)}) = 0$ for $k = 1,..., N-1$.

The above shows that $\left(\bigotimes_{i\in I}^{\omega_{M}} X_i\right) \cap \left(\sum_{k=1}^{M-1} \bigotimes_{i\in I}^{\omega_k} X_i\right) = \{0\}$ whenever $\omega_1,...,\omega_M$ are distinct elements in $\Omega_{I;X}$. 
On the other hand, for every $x\in \PI  X_i^\times$, one has $\Theta_X(x)\in \bigotimes_{i\in I}^{[x]_\sim} X_i$. 
These give the required equality. 
\end{prf}

\medskip

For any $F\in \KF$ and $u\in \PI  X_i^\times$, one has a linear map $$J_F^u\ :\ {\bigotimes}_{i\in F} X_i\ \longrightarrow\ \BOI^u X_i$$ 
given by $J_F^u(\otimes_{i\in F}\ \! x_i) := \otimes_{j\in I}\ \! \ti x_j$ ($x_i\in X_i$), where 
$\ti x_j := x_j$ when $j\in F$, and
$\ti x_j := u_j$ when $j\in I\setminus F$.

\medskip

For any $F,G\in \mathfrak{F}$ with $F\subseteq G$, a similar construction gives a linear map $J_{G;F}^u: \bigotimes_{i\in F} X_i \to \bigotimes_{i\in G} X_i$. 
It is clear that $\left(\bigotimes_{i\in F} X_i, J^u_{G;F}\right)_{F\subseteq G\in \mathfrak{F}}$ is an inductive system in the category of vector spaces with linear maps as morphisms. 

\medskip

\begin{prop}\label{prop:ju-inj}
(a) $J_F^u$ is injective for any $u\in \PI  X_i^\times$ and $F\in \mathfrak{F}$. 
Consequently, $\Theta_X(u)\neq 0$.

\smnoind
(b) The inductive limit of $\left(\bigotimes_{i\in F} X_i, J^u_{G;F}\right)_{F\subseteq G\in \mathfrak{F}}$ is $\left( \bigotimes_{i\in I}^u X_i, \{J_F^u\}_{F\in \mathfrak{F}}\right)$.
\end{prop}
\begin{prf}
(a) Suppose that $a\in \ker J_F^u$ and $\psi\in (\bigotimes_{i\in F} X_i)^*$. 
For each $j\in I\setminus F$, choose $f_j\in X_j^*$ with $f_j(u_j) = 1$. 
Remark \ref{rem:alt-con-ten-prod}(b) and the universal property give a linear map $\check \psi: \bigotimes_{i\in I} X_i \to \BC^{{\otimes I}}$ satisfying 
$$\check\psi(\otimes_{i\in I} x_i) = \psi(\otimes_{i\in F}\ \!x_i)\left(\otimes_{j\in I\setminus F}\ \! f_j(x_j)\right) \qquad (x\in \PI  X_i).$$ 
Thus, $\psi(a)(\tsi 1) = \check\psi(J_F^u(a)) = 0$, which implies that $a = 0$ (as $\psi$ is arbitrary) as required. 
On the other hand, if $i_0\in I$, then $\Theta_X(u) = J_{\{i_0\}}^u(u_{i_0})\neq 0$. 

\smnoind
(b) This follows directly from part (a). 
\end{prf}

\medskip

Part (b) of the above implies that $\Theta_X(C^\omega)$ is a basis for $\BOI^\omega X_i$, where $C^\omega$ is as defined in the following result. 

\medskip

\begin{cor}\label{cor:sub-sp}
(a) Let $c: \Omega_{I;X} \to \PI  X_i^\times$ be a cross-section.
For each $\omega\in \Omega_{I;X}$ and $\ii$, we pick a basis $B^\omega_i$ of $X_i$ that contains $c(\omega)_i$ and set
$$C^\omega := \{x\in \PI^{\omega}  X_i: x_i\in B^\omega_i, \forall i\in I\}.$$ 
If $C := \bigcup_{\omega\in \Omega_{I;X}} C^\omega$, 
then $\Theta_X(C)$ is a basis for $\bigotimes_{i\in I} X_i$.

\smnoind
(b) If $\Phi_i:X_i\to Y_i$ is an injective linear map ($i\in I$),  
the induced linear map $\BOI \Phi_i: \BOI X_i \to \BOI Y_i$ is injective. 
\end{cor}

\medskip

\begin{prop}\label{prop:inj-TBOIM}
The map $\Psi: \BOI L(X_i;Y_i) \to L(\BOI X_i; \BOI Y_i)$ (given by the universal property) is injective. 
\end{prop}
\begin{prf}
Suppose that $T^{(1)}, ..., T^{(n)}\in \PI L(X_i;Y_i)^\times$ are mutually inequivalent elements (under $\sim$), $F\in \mathfrak{F}$, $R^{(1)},...,R^{(n)}\in \bigotimes_{i\in F} L(X_i;Y_i)$ with $S^{(k)} := J_F^{T^{(k)}}(R^{(k)})$ ($k=1,...,n$) satisfying 
\begin{equation*}
\Psi\big({\sum}_{k=1}^n S^{(k)}\big)\ =\ 0.
\end{equation*}
Using an induction argument, it suffices to show that $S^{(1)} = 0$. 

If $n=1$, we take any $x\in \PI X_i^\times$ with $T_i^{(1)}x_i\neq 0$ ($\ii$). 
If $n >1$, we claim that there is $x\in \PI X_i^\times$ such that $$[T^{(1)}_ix_i]_\ii\in \PI Y_i^\times\ \text{ and }\ [T^{(k)}_ix_i]_\ii \nsim [T^{(1)}_ix_i]_\ii \ \ (k =2,...,n).$$ 
In fact, let $I^k:= \{i\in I: T^{(k)}_i \neq T^{(1)}_i\}$, which is an infinite set for any $k=2,...,n$. 
For any $i\in I$, we put $N_i := \{k\in \{2,,..,n\}: i\in I^k\}$ and  
take any $x_i \in X_i \setminus \big(\bigcup_{k\in N_i} \ker (T_i^{(k)} - T_i^{(1)})\cup\ker T_i^{(1)}\big)$ (note that $X_i$ cannot be a finite union of proper subspaces).  
Thus, $T^{(1)}_ix_i \neq 0$ (for each $\ii$) and $T^{(k)}_ix_i \neq T^{(1)}_ix_i$ (for $k\in\{2,..,n\}$ and $i\in I^k$). 

Now, we have 
$$\Psi(S^{(1)})\big(\BOI^x X_i\big) \cap \Big({\sum}_{k=2}^n \Psi(S^{(k)})\big(\BOI^x X_i\big)\Big) = (0)$$ 
by Theorem \ref{thm:bas-inf-ten} and the fact that $\Psi(S^{(l)})\big(\BOI^x X_i\big)\in \BOI^{y^{(l)}} Y_i$, 
where $y^{(l)}_i = T^{(l)}_ix_i$ ($\ii; l= 1,...,n$).
Consequently, $\Psi(S^{(1)})\big|_{\BOI^x X_i} = 0$.
As $T^{(1)}_ix_i\neq 0$ ($i\in I$), it is easy to see that $R^{(1)} = 0$ as required. 
\end{prf}

\medskip

Note that $\Psi$ is not surjective even if $X_i=Y_i =\BC$ ($\ii$) since in this case, $\Psi$ is a homomorphism and $\BOI \BC$ is commutative while $L(\BOI \BC)$ is not. 

\medskip

The following result follows from Proposition \ref{prop:exist-ten-prod}(c), Corollary \ref{cor:sub-sp}(b) and Proposition \ref{prop:inj-TBOIM}, which say that an infinite tensor product of vector spaces is automatically a faithful module over a big commutative algebra. 

\medskip

\begin{cor}\label{cor:ten-prod-as-mod}
If $X_i$ is a faithful $A_i$-module ($\ii$), then $\BOI X_i$ is a faithful $\BOI A_i$-module. 
In particular, $\bigotimes_{i\in I} Y_i$ is a faithful unital $\BC^{\otimes I}$-module. 
\end{cor}

\medskip

\begin{eg}\label{eg:u-ten-prod}
(a) If $\beta \in \PI \mathbb{C}^\times$, then ${\bigotimes}_{i\in I}^\beta \mathbb{C} = \mathbb{C}\cdot \tsi \beta_i$.  
In fact, for any $F\in \mathfrak{F}$ and $\mu_i\in \mathbb{C}$ ($i\in F$), we have 
$J_F^\beta({\otimes}_{i\in F}\ \!\mu_i) = \left({\Pi}_{i\in F} \ \!\mu_i/\beta_i\right)(\tsi \beta_i)$. 

\smnoind
(b) Let $n\in \BN$, $I_1, ..., I_n$ be infinite disjoint subsets of $I$ with $I = \bigcup_{k=1}^n I_k$ and $\overline{\beta}=(\beta_1, ..., \beta_n)\in (\BC^\times)^n$. 
Define $\widetilde{\beta}\in \PI \BC^\times$ by $\widetilde{\beta}_i = \beta_k$ whenever $i\in I_k$. 
Then $\overline{\beta} \mapsto [\widetilde{\beta}]_\sim$
is an injective group homomorphism from $(\BC^\times)^n$ to $\Omega_{I;\BC}$. 

\smnoind
(c) Let $G$ be a subgroup of $\BT^n\subseteq (\BC^\times)^n$ (where $\BT:=\{t\in \BC: \abs{t} =1\}$). 
If $\overline{\beta^{(1)}},...,\overline{\beta^{(m)}}$ are distinct elements in $G$ and $\widetilde{\beta^{(1)}},...,\widetilde{\beta^{(m)}}\in \PI \BC^\times$ are as in part (b), then $\tsi \widetilde{\beta^{(1)}_i},...,\tsi \widetilde{\beta^{(m)}_i}$ are linearly independent in $\BC^{\otimes I}$. 
Therefore, the $^*$-subalgebra of $\BC^{\otimes I}$ generated by $\{\tsi \widetilde{\beta}_i: \overline{\beta}\in G\}$ is $^*$-isomorphic to 
the group algebra $\BC[G]$. 
\end{eg}

\medskip

As $\tsi \alpha_i = (\PI \alpha_i) (\tsi 1)$ if $\alpha_i =1$ e.f., one may regard $\tsi \alpha_i$ as a generalisation of the product. 
In this case, one can consider infinite products like $(-1)^I$. 

\medskip
\section{Tensor products of unital $^*$-algebras}

\medskip

\emph{Throughout this section, $A_i$ is a unital $^*$-algebra with identity $e_i$ ($\ii$), and we set $\Omega^\ut_{I;A}:=\PI   U_{A_i}/\sim$.}

\medskip

Notice that in this case, $\Omega_{I;A}$ is a $^*$-semi-group with identity and $\Omega^\ut_{I;A}$ can be regarded as a subgroup of $\Omega_{I;A}$ with the inverse being the involution on $\Omega_{I;A}$. 
Moreover, $\bigotimes_{i\in I} A_i$ is a $\Omega_{I;A}$-graded $^*$-algebra in the sense that for any $\omega, \omega'\in \Omega_{I;A}$, 
\begin{equation}\label{eqt:grading}
\Big({\bigotimes}_{i\in I}^{\omega} A_i\Big) \cdot \Big({\bigotimes}_{i\in I}^{\omega'} A_i\Big)
\ \subseteq\ {\bigotimes}_{i\in I}^{\omega\omega'} A_i 
\ \text{ and }\  
\Big({\bigotimes}_{i\in I}^{\omega} A_i\Big)^*
\ \subseteq\ {\bigotimes}_{i\in I}^{\omega^*} A_i. 
\end{equation}

\medskip

By Proposition \ref{prop:ju-inj}(b), $\BOI^e A_i$ can be identified with the unital $^*$-algebra inductive limit of finite tensor products of $A_i$. 
We will study the following $^*$-subalgebra that contains $\BOI^e A_i$: 
$$\BOI^\ut A_i\ :=\ {\bigoplus}_{\omega\in \Omega^\ut_{I;A}} \BOI^\omega A_i.$$ 
The motivation behind the consideration of this subalgebra is partially from Example \ref{eg:gp-alg-inf-prod-gp}(a) below, and partially because it has good representations (see the discussion after Proposition \ref{prop:ten-prod-st-alg} below). 
Moreover, if all $A_i$ are linear spans of $U_{A_i}$ (in particular, if they are $C^*$-algebras or group algebras), then $\BOI^\ut A_i$ is the linear span of $\Theta_A(\PI  U_{A_i})$. 
If $A_i = A$ for all $i\in I$, we denote $A^{\otimes I}_\ut := \BOI^\ut A_i$. 

\medskip

\begin{eg}\label{eg:gp-alg-inf-prod-gp}
(a) Let $G_i$ be a group and $\BC[G_i]$ be its group algebra ($i\in I$). 
If $\Lambda: \PI  G_i \to \PI  U_{\BC[G_i]}$ is the canonical map, then $\lambda := \Theta_{\BC[G]} \circ \Lambda$ gives a $^*$-isomorphism from $\BC[\PI  G_i]$ to the $^*$-subalgebra 
$$\BOI^{\Lambda(\PI G_i)} \BC[G_i]\ := \ {\sum}_{t\in \PI G_i} \BOI^{\Lambda(t)} \BC[G_i] \ \subseteq \ \BOI^\ut \BC[G_i].$$ 
In fact, $\lambda$ induces a $^*$-homomorphism from $\BC[\PI  G_i]$ to $\BOI^\ut \BC[G_i]$. 
Let $q: {\Pi}_\ii G_i \to {\Pi}_\ii G_i/{\oplus}_{i\in I} G_i$ be the quotient map. 
For a fixed $s\in \PI  G_i$, if we set 
$${\bigoplus}_{i\in I}^s G_i\ :=\ \big\{t\in \PI  G_i: q(t) =  q(s)\big\},$$ 
then $s^{-1} \big(\bigoplus_{i\in I}^s G_i \big) = \bigoplus_{i\in I} G_i$.
Thus, $\{\lambda(t): t\in \bigoplus_{i\in I}^s G_i\}$ is a set of linearly independent elements in $\BOI \BC[G_i]$ (as $\lambda|_{\BC[\bigoplus_{i\in I} G_i]}$ is a bijection onto $\BOI^e \BC[G_i]$). 
On the other hand, if $s^{(1)}, ..., s^{(n)}\in \PI  G_i$ such that $q(s^{(k)})\neq q(s^{(l)})$ whenever $k\neq l$, then $\lambda(s^{(1)}),...,\lambda(s^{(n)})$ are linearly independent in $\BOI\BC[G_i]$ (see Theorem \ref{thm:bas-inf-ten}). 
Consequently, $\{\lambda(t): t\in \PI G_i\}$ form a basis for $\BOI^{\Lambda(\PI G_i)} \BC[G_i]$. 

\smnoind
(b) It is well-known that there is a twisted action $(\alpha, u)$, in the sense of Busby and Smith, of $\Omega_{I;G}:= {\Pi}_\ii G_i/{\oplus}_{i\in I} G_i$ on $\BC[\bigoplus_{i\in I} G_i]\cong \BOI^e \BC[G_i]$ (see \cite[2.1]{BS}) such that 
$\BC[\PI G_i]$ is $^*$-isomorphic to the algebraic crossed-product $\BOI^e \BC[G_i]\rtimes_{\alpha,u} \Omega_{I;G}$. 
\end{eg}

\medskip

There is a canonical action $\Xi$ of $\PI  U_{A_i}$ on $\BOI^\ut A_i$ given by inner-automorphisms, i.e. 
$$\Xi_u(a) := (\tsi u_i)\cdot a\cdot (\tsi u^*_i) 
\qquad \big(u\in \PI  U_{A_i}; a\in \BOI^\ut A_i\big).$$ 
This induces an action $\Xi^e$ of $\PI  U_{A_i}$ on the subalgebra $\BOI^e A_i$. 
The following result gives an identification of $\BOI^\ut A_i$ as the algebraic crossed-product (see e.g. \cite[p.166]{RSW}) of a cocycle twisted action (i.e. a twisted action in the sense of Busby and Smith) of $\Omega^\ut_{I;A}$ on $\BOI^e A_i$ induced by $\Xi^e$. 

\medskip

Before we give this result, let us recall that an abelian group $G$ is \emph{divisible} if for any $g\in G$ and $n\in \BN$, there is $h\in G$ with $g = h^n$. 

\medskip

\begin{thm}\label{thm:boiut-tw-cr-pd}
(a) There is a cocycle twisted action $(\check\Xi, m)$ of $\Omega^\ut_{I;A}$ on $\BOI^e A_i$ such that $\BOI^\ut A_i$ is $\Omega^\ut_{I;A}$-graded $^*$-isomorphic to $(\BOI^e A_i)\rtimes_{\check \Xi, m} \Omega^\ut_{I;A}$. 

\smnoind
(b) Suppose that all $A_i$ are commutative. 
If $\BOI^e A_i$ is a unital $^*$-subalgebra of a commutative $^*$-algebra $B$ with $U_B$ being divisible, $\BOI^\ut A_i$ is $\Omega_{I;A}^\ut$-graded $^*$-isomorphic to a unital $^*$-subalgebra of $B\otimes \BC[\Omega_{I;A}^\ut]$. 
If $U_{\BOI^e A_i}$ is itself divisible, $\BOI^\ut A_i \cong (\BOI^e A_i)\otimes \BC[\Omega_{I;A}^\ut]$ as $\Omega_{I;A}^\ut$-graded $^*$-algebras. 
\end{thm}
\begin{prf}
Let $c:\Omega^\ut_{I;A} \to \PI  U_{A_i}$ be a cross-section with $c([e]_\sim) = e$.

\smnoind
(a) For any $\mu,\nu\in \Omega^\ut_{I;A}$, we set
$$\check \Xi_\mu\ :=\ \Xi^e_{c(\mu)} \quad \text{and} \quad m(\mu,\nu)\ :=\ \tsi c(\mu)_ic(\nu)_ic(\mu\nu)^{-1}_i.$$
As $c(\mu)c(\nu)\sim c(\mu\nu)$, we have $m(\mu,\nu)\in \BOI^e A_i$. 
It is easy to check that $(\check \Xi, m)$ is a twisted action in the sense of Busby and Smith. 
Furthermore, we define $\Psi: (\BOI^e A_i)\rtimes_{\check \Xi, m} \Omega^\ut_{I;A} \to \BOI^\ut A_i$ by 
$$\Psi(f)\ :=\ {\sum}_{\omega\in \Omega^\ut_{I;A}} f(\omega)(\tsi c(\omega)_i)
\qquad \big(f\in (\BOI^e A_i)\rtimes_{\check \Xi, m} \Omega^\ut_{I;A}\big).
$$
It is not hard to verify that $\Psi$ is a bijective $\Omega^\ut_{I;A}$-graded $^*$-homomorphism. 

\smnoind
(b) Let $\PI^e U_{A_i} := \PI^e A_i \cap \PI U_{A_i}$. 
By the Baer's theorem, $\Theta_A|_{\PI^e U_{A_i}}$ can be extended to a group homomorphism $\varphi: \PI U_{A_i} \to U_B$.
Since
$$\varphi(c(\mu))\varphi(c(\nu))\varphi(c(\mu\nu))^{-1}
\ = \ \tsi c(\mu)_ic(\nu)_ic(\mu\nu)^{-1}_i
\qquad (\mu,\nu\in \Omega^\ut_{I;A}),$$
the map $\Phi: \BOI^\ut A_i \to B \otimes \BC[\Omega^\ut_{I;A}]$ given by 
\begin{equation}\label{eqt:Phi}
\Phi(a)
 := (a \cdot \tsi c(\omega)^{-1}_i)\varphi(c(\omega)) \otimes \lambda(\omega)
\qquad
\big(a\in \BOI^\omega A_i; \omega\in \Omega^\ut_{I;A}\big)
\end{equation}
is a $\Omega^\ut_{I;A}$-graded $^*$-homomorphism. 
If ${\sum}_{\omega\in \Omega^\ut_{I;A}} a^\omega\in \ker \Phi$ (with $a^\omega\in \BOI^\omega A_i$), then for every $\omega\in \Omega^\ut_{I;A}$, one has $(a^\omega \cdot \tsi c(\omega)^{-1}_i)\varphi(c(\omega)) = 0$, which implies $a^\omega = 0$, and hence $\Phi$ is injective. 
The image of $\Phi$ is the linear span of 
$$\big\{b\varphi(c(\omega))\otimes \lambda(\omega): b\in \BOI^eA_i; \omega\in \Omega^\ut_{I;A}\big\},$$
and it is clear that $\Phi$ is surjective if $B = \BOI^eA_i$. 
\end{prf}

\medskip

\begin{rem}\label{rem:ut-prod}
(a) The cocycle twisted action $(\check\Xi, m)$ depend on the choice of a cross-section, and different cross-sections may give different twisted actions (although their crossed-products are all isomorphic). 
On the other hand, the map $\Phi$ in part (b) also depends on the choice of a cross-section as well as the choice of an extension of $\Theta_A|_{\PI^e U_{A_i}}$. 

\smnoind
(b) If $S_i$ is a set and $A_i$ is a $^*$-subalgebra of $\ell^\infty(S_i)$ ($i\in I$), then by Theorem \ref{thm:boiut-tw-cr-pd}(b), $\BOI^\ut A_i$ is a $^*$-subalgebra of $\ell^\infty(\PI S_i)\otimes \BC[\Omega^\ut_{I;A}]$. 
Our first proof for this fact use \cite[18.4]{Cal} and \cite[7.1]{EM}.  

\smnoind
(c) If all $A_i$ are commutative, then 
$\BOI^\ut A_i \cong (\BOI^e A_i)\otimes \BC[\Omega_{I;A}^\ut]$ as $\Omega_{I;A}^\ut$-graded $^*$-algebras if and only if there is a group homomorphism $\pi: \Omega^\ut_{I;A} \to U_{\BOI^\ut A_i}$ such that $\pi(\omega)\in \BOI^\omega A_i$ ($\omega\in \Omega^\ut_{I;A}$). 
In fact, if such a $\pi$ exists, one may replace $(a\cdot \tsi c(\omega)^{-1}_i)\varphi(c(\omega))$ in \eqref{eqt:Phi} with $a\pi(\omega^{-1})$ and show that the corresponding $\Phi$ is a $^*$-isomorphism. 
\end{rem}

\medskip

Clearly, the second statement of Theorem \ref{thm:boiut-tw-cr-pd}(b) applies to the case when $A_i = \BC^{n_i}$ for some $n_i\in \BN$ ($\ii$). 
In particular, Theorem \ref{thm:boiut-tw-cr-pd}(b) and its argument give the following corollary. 

\medskip

\begin{cor}\label{cor:boiut-C}
If $\varphi_1$ is as in Example \ref{eg:multi-linear}(a) and $\varphi: \PI \BT\to \BT$ is a group homomorphism that extends $\varphi_1|_{\PI^1 \BT}$ (it existence is guaranteed by the Baer's theorem), then $\Phi(\tsi \alpha_i) := \varphi(\alpha) \lambda([\alpha]_\sim)$ ($\alpha \in \PI \BT$) is a well-defined $^*$-isomorphism from 
$\BC^{\otimes I}_\ut$ onto $\BC[\Omega^\ut_{I;\BC}]$.
\end{cor}

\medskip

Conversely, it is clear that if $\varphi: \PI \BT\to \BT$ is any map such that $\Phi$ as defined in the above is a well-defined $^*$-isomorphism, then $\varphi$ is a group homomorphism extending  $\varphi_1|_{\PI^1 \BT}$. 
On the other hand, there is a simpler proof for Corollary \ref{cor:boiut-C}. 
In fact, for $\alpha,\beta\in \PI \BT$ with $\alpha\sim \beta$, one has $\varphi(\alpha)^{-1}\cdot \tsi \alpha_i = \varphi(\beta)^{-1}\cdot \tsi \beta_i$. 
Thus, $[\alpha]_\sim \mapsto \varphi(\alpha)^{-1} \cdot \tsi \alpha_i$ is a well-defined group homomorphism from $\Omega^\ut_{I;\BC}$ to $U_{\BC^{\otimes I}_\ut}$ such that $\{\varphi(\alpha)^{-1}\cdot \tsi \alpha_i: [\alpha]_\sim \in \Omega^\ut_{I;\BC}\}$ is a basis for $\BC^{\otimes I}_\ut$. 

\medskip

\begin{eg}\label{eg:mT-Omega}
For any subgroup $G\subseteq \mathbb{T}^n$, the algebra in Example \ref{eg:u-ten-prod}(c) is a $^*$-subalgebras of $\BC^{\otimes I}_\ut$. 
\end{eg}

\medskip

In the remainder of this section, we will show that the center of $\BOI^\ut A_i$ is the tensor product of centers of $A_i$ when $A_i = {\rm span}\ \!U_{A_i}$ for all $i\in I$. 

\medskip

If $A$ is an algebra and $G$ is a group, we denote by $Z(A)$ and $Z(G)$ the center of $A$ and the center of $G$ respectively. 
Clearly, the inclusion $\PI  U_{Z(A_i)} \subseteq \PI  U_{A_i}$ induces an injective group homomorphism from $\Omega^\ut_{I;Z(A)}$ to $\Omega^\ut_{I;A}$ and we regard the former as a subgroup of the later. 

\medskip

\begin{thm}\label{thm:center}
Suppose that there is $F_0\in \KF$ with $A_i = {\rm span}\ \!U_{A_i}$ for any $i\in I_0:= I\setminus F_0$. 

\smnoind
(a) $Z(\Omega^\ut_{I;A}) = \Omega^\ut_{I;Z(A)}$. 
Moreover, $Z(\Omega^\ut_{I;A}) = \Omega^\ut_{I;A}$ if and only if all but a finite number of $A_i$ are commutative. 

\smnoind
(b) Every element in $\Omega^\ut_{I;A}\setminus Z(\Omega^\ut_{I;A})$ has an infinite conjugacy class. 

\smnoind
(c) $Z\big(\BOI^\ut A_i\big) = \BOI^\ut Z(A_i)$.
\end{thm}
\begin{prf}
(a) It is obvious that $\Omega^\ut_{I;Z(A)}\subseteq Z(\Omega^\ut_{I;A})$. 
Suppose $u\in \PI  U_{A_i}$ with $[u]_\sim \notin \Omega^\ut_{I;Z(A)}$. 
There is an infinite subset $J\subseteq I_0$ such that $u_i\notin Z(A_i)$ ($i\in J$). 
For each $i\in J$, one can find $v_i\in U_{A_i}$ such that $u_iv_i \neq v_iu_i$. 
For any $i\in I\setminus J$, we put $v_i = e_{i}$. 
Then $[v]_\sim \in \Omega^\ut_{I;A}$ and $[u]_\sim [v]_\sim \neq [v]_\sim[u]_\sim$. 
Consequently, $[u]_\sim \notin Z(\Omega^\ut_{I;A})$. 
This argument also shows that if the set $\{i\in I: Z(A_i)\neq A_i\}$ is infinite, then $Z(\Omega^\ut_{I;A}) \neq \Omega^\ut_{I;A}$.
Conversely, it is clear that $\Omega^\ut_{I;Z(A)} = \Omega^\ut_{I;A}$ if all but a finite numbers of $A_i$ are commutative. 

\smnoind
(b) Suppose that $[u]_\sim \in \Omega^\ut_{I;A}\setminus Z(\Omega^\ut_{I;A})$ and $\{i_n\}_{n\in \BN}$ is a sequence of distinct elements in $I_0$ such that $u_{i_n}\notin Z(A_{i_n})$ ($n\in \BN$). 
For each $n\in \BN$, choose $v_{i_n}\in U_{A_{i_n}}$ with $v_{i_n}u_{i_n}v_{i_n}^* \neq u_{i_n}$. 
For any prime number $p$, we 
set $w^{(p)}_{i_n} := v_{i_n}$ ($n\in \BN p$), and $w^{(p)}_i := e_{i}$ if $i\in I\setminus \{i_{n}: n\in \BN p\}$. 
If $p$ and $q$ are distinct prime numbers, then
$$w^{(q)}_{i_n} u_{i_n} (w^{(q)}_{i_n})^*
\ =\ u_{i_n} 
\ \neq\ w^{(p)}_{i_n} u_{i_n} (w^{(p)}_{i_n})^*
\qquad (n\in \BN p \setminus \BN q).$$
Consequently, $w^{(q)} u (w^{(q)})^*
\nsim w^{(p)} u (w^{(p)})^*$, and the conjugacy class of $[u]_\sim$ is infinite. 

\smnoind
(c) Since $Z(\BOI^\ut A_i) = {\bigotimes}_{i\in F_0} Z(A_i)\otimes Z({\bigotimes}_{i\in I_0}^\ut A_i)$, we may assume that $A_i = {\rm span}\ \! U_{A_i}$ for all $\ii$. 
In this case, $Z(\BOI^\ut A_i) = \big(\BOI^\ut A_i\big)^\Xi$, where 
$\big(\BOI^\ut A_i\big)^\Xi$ is the fixed point algebra of the action $\Xi$ as defined above. 
Moreover, one has $\BOI^\ut Z(A_i) \subseteq Z(\BOI^\ut A_i)$ and it remains to show that $\big(\BOI^\ut A_i\big)^\Xi \subseteq \BOI^\ut Z(A_i)$. 

Let $v^{(1)},...,v^{(n)}\in \PI  U_{A_i}$ be mutually inequivalent elements, $F\in \mathfrak{F}$ and $b_1,...,b_n\in \BOF A_i\setminus \{0\}$ such that $a := {\sum}_{k=1}^n J_F^{v^{(k)}}(b_k)\in \big(\BOI^\ut A_i\big)^\Xi$. 
We first claim that $[v^{(k)}]_\sim \in \Omega^\ut_{I;Z(A)}$ ($k = 1,...,n$). 
Suppose on the contrary that $[v^{(1)}]_\sim \notin \Omega^\ut_{I;Z(A)} = Z(\Omega^\ut_{I;A})$. 
For every $u\in \PI  U_{A_i}$, one has 
$$\Xi_u\big(J_F^{v^{(1)}}(b_k)\big)\ \in \ \big(\BOI^{[uv^{(1)}u^*]_\sim} A_i\big) \setminus \{0\}.$$ 
As $\Xi_u(a) = a$, we see that $[uv^{(1)}u^*]_\sim \in \{[v^{(1)}]_\sim, ..., [v^{(n)}]_\sim\}$, which contradicts the fact that $\{[uv^{(1)}u^*]_\sim: [u]_\sim \in \Omega^\ut_{I;A}\}$ is an infinite set (by part (b)). 

By enlarging $F$, we may assume that $v^{(k)}\in \PI  U_{Z(A_i)}$ ($k = 1,...,n$).
For each $u\in \PI  U_{A_i}$ and $k\in \{1,...,n\}$, one has $\Xi_u(J_F^{v^{(k)}}(b_k)) = J_F^{v^{(k)}}(b_k)$ and so, $b_k \in Z(\BOF A_i)$. 
Therefore, $a \in \BOI^\ut Z(A_i)$ as expected. 
\end{prf}

\medskip

The readers should notice that $\BOI^\ut Z(A_i)$ equals $\bigoplus_{\omega \in Z(\Omega^\ut_{I;A})} \BOI^\omega Z(A_i)$ instead of $\bigoplus_{\omega \in \Omega^\ut_{I;A}} \BOI^\omega Z(A_i)$ (strictly speaking, the later object does not make sense). 

\medskip

\begin{eg}
(a) If $n_i\in \BN$ ($i\in I$), then 
$Z\big(\BOI^\ut M_{n_i}(\BC)\big) \cong \BC^{\otimes I}_\ut$. 

\smnoind
(b) If $G_i$ are icc groups, then $Z(\BOI^\ut \BC[G_i]) \cong \BC^{\otimes I}_\ut$ canonically. 
\end{eg}

\medskip

We end this section with the following brief discussion on the non-unital case. 
Suppose that $\{A_i\}_\ii$ is a family of $^*$-algebras, not necessarily unital. 
If $M(A_i)$ is the double centraliser algebra of $A_i$ ($\ii$), we define an ideal, $\BOI^\ut A_i$, of $\BOI^\ut M(A_i)$ as follows: 
$$\BOI^\ut A_i
\ :=\ {\rm span}\ \!\big\{J_F^u(a): F\in \KF; a\in {\bigotimes}_{i\in F} A_i; u\in \PI U_{M(A_i)}\big\}.$$ 
In general, $\BOI^\ut A_i$ is not a subset of $\BOI A_i$. 
In a similar fashion, we define 
$$\BOI^e A_i\ :=\ {\rm span}\ \!\big\{J_F^u(a): F\in \KF; a\in {\bigotimes}_{i\in F} A_i; u\in \PI U_{M(A_i)}; u\sim e\big\},$$ which is an ideal of $\BOI^e M(A_i)$. 
By the proof of Theorem \ref{thm:boiut-tw-cr-pd}(a), one may identify $\BOI^\ut A_i$ as the ideal of $(\BOI^e M(A_i))\rtimes_{\check \Xi, m} \Omega^\ut_{I;M(A)}$ consisting of functions from $\Omega^\ut_{I;M(A)}$ to $\BOI^e A_i$ having finite supports. 

\medskip

\section{Tensor products of inner-product spaces}

\medskip

\emph{Throughout this section, $(H_i, \langle \cdot, \cdot\rangle)$ is a non-zero inner-product space ($i\in I$). 
Moreover, we denote $\Omega^\un_{I;H} := \PI  \sph(H_i)/\sim$. }

\medskip

If $B$ is a unital $^*$-algebra and $X$ is a unital left $B$-module, a map $\langle \cdot, \cdot \rangle_B : X\times X \to B$ is called a \emph{(left) Hermitian $B$-form on $X$} if 
$\langle ax + y, z \rangle_B = a \langle x, z \rangle_B + \langle y, z \rangle_B$ and $\langle x, y \rangle_B^* = \langle y, x \rangle_B$ ($x,y,z\in X; a\in B$).
It is easy to see that a Hermitian $B$-form on $X$ can be regarded as a $B$-bimodule map $\theta : X\otimes \ti X \to B$ satisfying $\theta(x\otimes \ti y)^* = \theta (y\otimes \ti x)$ (where $\ti X$ is the conjugate vector spaces of $X$ regarded as a unital right $B$-module in the canonical way). 
Consequently, part (a) of the following result follows readily from the universal property of tensor products, while part (b) is easily verified. 

\medskip

\begin{prop}\label{prop:ten-inn-prod}
(a) There is a Hermitian $\BC^{\otimes I}$-form on $\bigotimes_{i\in I} H_i$ such that 
$\left< \tsi x_i, \tsi y_i \right>_{\BC^{\otimes I}}
:= \otimes_{i\in I} \left< x_i, y_i \right>$ ($x,y\in \PI  H_i$).

\smnoind
(b) For a fixed $\mu\in \Omega^\un_{I;H}$, one has $\left< \Theta_H(x), \Theta_H(y)\right>_{\BC^{\otimes I}}
= \PI  \langle x_i, y_i \rangle (\tsi 1)$ $(x,y\in \PI^\mu  H_i)$. 
This induces an inner-product on $\BOI^\mu H_i$ which coincides with the one given by the inductive limit of $\big(\bigotimes_{i\in F} H_i, J^\mu_{G;F}\big)_{F\subseteq G\in \mathfrak{F}}$, in the category of inner-product spaces with isometries as morphisms. 
\end{prop}

\medskip

We want to construct a nice inner-product space from the above Hermitian $\BC^{\otimes I}$-form. 
A naive thought is to appeal to a construction in Hilbert $C^*$-modules that produces a Hilbert space from a positive linear functional on $\BC^{\otimes I}$. 
However, the difficulty is that there is no canonical order structure on $\BC^{\otimes I}$. 
Nevertheless, we will do a similar construction using the functional $\phi_1$ in Example \ref{eg:multi-linear}(a). 
In this case, one can only consider a subspace of $\BOI H_i$ (see Example \ref{eg:bad-expan} below). 

\medskip

\begin{lem}\label{lem:phi0-innpro}
Suppose that $\bigotimes_{i\in I}^\ct H_i:= {\rm span}\ \!\Theta_H(\PI  B_1(H_i))$, $\bigotimes_{i\in I}^\un H_i:= {\rm span}\ \!\Theta_H(\PI  \sph(H_i))$ and  
$$\langle \xi, \eta \rangle_{\phi_1}
\ :=\ \phi_1( \langle \xi, \eta \rangle_{\BC^{\otimes I}})
\qquad \big(\xi,\eta\in \BOI H_i\big).$$ 

\smnoind
(a) For any $\mu\in \Omega^\un_{I;H}$, the restriction of 
$\langle \cdot, \cdot \rangle_{\phi_1}$ on $\bigotimes_{i\in I}^\mu H_i\times \bigotimes_{i\in I}^\mu H_i$ coincides with the inner-product in Proposition \ref{prop:ten-inn-prod}(b). 

\smnoind
(b) $\left< \cdot, \cdot \right>_{\phi_1}$ is a positive sesquilinear form on $\bigotimes_{i\in I}^\ct H_i$ and is an 
inner-product on $\bigotimes_{i\in I}^\un H_i$.
Moreover, if 
\begin{equation*}
K\ :=\ \Big\{y\in \BOI^\ct H_i: \langle x, y\rangle_{\phi_1} = 0, \forall x\in \BOI^\ct H_i \Big\},
\end{equation*}
then $\bigotimes_{i\in I}^\ct H_i = K \oplus \bigotimes_{i\in I}^\un H_i$ (as vector spaces). 

\smnoind
(c) If $I = I_1 \cup I_2$ and $I_1 \cap I_2 = \emptyset$, then $\bigotimes_{i\in I}^\un H_i = (\bigotimes_{i\in I_1}^\un H_i) \otimes (\bigotimes_{j\in I_2}^\un H_j)$ as inner-product spaces. 
\end{lem}
\begin{prf}
(a) This part is clear. 

\smnoind
(b) It is obvious that $\left< \cdot, \cdot \right>_{\phi_1}$ is a sesquilinear form on $\bigotimes_{i\in I}^\ct H_i$. 
Let 
$$E\ :=\ \big\{x\in \PI  B_1(H_i): \|x_i\| < 1 \text{ for an infinite number of } i\in I\big\}$$ 
and $\ti K:= {\rm span}\ \! \Theta_H(E)$. 
Clearly, $\bigotimes_{i\in I}^\ct H_i = \ti K\oplus \bigotimes_{i\in I}^\un H_i$.
Moreover, if $u\in \PI B_1(H_i)$ and $v\in E$, then $\langle u_i, v_i \rangle \neq 1$ for an infinite number of $i\in I$, which implies that $\langle \tsi u_i, \tsi v_i\rangle_{\phi_1} = 0$. 
Consequently, $\ti K\subseteq K$. 

We claim that $\left< \xi, \xi \right>_{\phi_1} \geq 0$ for any $\xi\in \bigotimes_{i\in I}^\ct H_i$. 
Suppose that $\xi = \sum_{k=1}^n \lambda_k \tsi u^{(k)}_i$ with $\lambda_1,...,\lambda_n\in \BC$ and $u^{(1)}, ..., u^{(n)}\in \PI  B_1(H_i)$. 
Then 
$$\left< \xi, \xi \right>_{\phi_1}
\ = \ {\sum}_{k,l=1}^n \lambda_k \bar\lambda_l \phi_1\big(\otimes_{i\in I} \langle u^{(k)}_i, u^{(l)}_i\rangle\big).$$
As in the above, $\phi_1\big(\otimes_{i\in I} \langle u^{(k)}_i, u^{(l)}_i\rangle\big) = 0$ if either $u^{(k)}$  or $u^{(l)}$ is in $E$. 
Thus, by rescaling, we may assume that 
$$u^{(1)}, ..., u^{(n)}\in \PI  \sph(H_i).$$
Furthermore, we assume that there exist $0 = n_0  < \cdots < n_m = n$ such that 
$u^{(n_p+1)} \sim \cdots \sim u^{(n_{p+1})}$ for all $p \in  \{0,...,m-1\}$, but $u^{(n_p)}\nsim u^{(n_q)}$ whenever $1\leq p\neq q\leq m$. 
It is not hard to check that $u^{(k)}\sim u^{(l)}$ if and only if $\langle u^{(k)}_i, u^{(l)}_i \rangle = 1$ e.f. (as $\|u^{(k)}_i\|, \|u^{(l)}_i\|\leq 1$).
Consequently, if $1\leq p\neq q\leq m$, \begin{equation}\label{eqt:phi0-disj}
\phi_1\big(\otimes_{i\in I} \langle u^{(k)}_i, u^{(l)}_i\rangle \big) = 0
\quad \text{when } n_p < k \leq n_{p+1} \text{ and } n_q < l \leq n_{q+1}.
\end{equation}
Therefore, in order to show $\left< \xi, \xi \right>_{\phi_1} \geq 0$, it suffices to consider the case when $u^{(k)} \sim u^{(l)}$ for all $k,l \in \{1,...,n\}$, which is the same as $\xi\in \bigotimes_{i\in I}^{u^{(1)}} H_i$. 
Thus, $\left< \xi, \xi \right>_{\phi_1} \geq 0$ by part (a). 

Next, we show that $\left< \cdot, \cdot \right>_{\phi_1}$ is an inner-product on $\bigotimes_{i\in I}^\un H_i$. 
Suppose that $\xi = \sum_{k=1}^n \lambda_k \tsi u^{(k)}_i$ with $\lambda_1,...,\lambda_n\in \BC$ and $u^{(1)}, ..., u^{(n)}\in \PI  \sph(H_i)$ such that $\left< \xi, \xi \right>_{\phi_1} = 0$. 
If $n_0,...,n_m$ are as in the above, then  
$$\phi_1\Big( \big< {\sum}_{k=n_{p}+1}^{n_{p+1}} \lambda_k \tsi u^{(k)}_i, {\sum}_{l=n_{q}+1}^{n_{q+1}} \lambda_l \tsi u^{(l)}_i\big>_{\BC^{\otimes I}} \Big) 
\ = \ 0,$$
because of \eqref{eqt:phi0-disj} and the positivity of $\langle \cdot, \cdot\rangle_{\phi_1}$.
Hence, we may assume $u^{(k)} \sim u^{(l)}$ for all $k,l \in \{1,...,n\}$, and apply part (a) to conclude that $\xi = 0$. 

Finally, as $\left< \cdot, \cdot \right>_{\phi_1}$ is an inner-product on $\BOI^\un H_i$ and we have both $\bigotimes_{i\in I}^\ct H_i = \ti K\oplus \bigotimes_{i\in I}^\un H_i$ and $\ti K \subseteq K$, we obtain  
$K \subseteq \ti K$ as well. 

\smnoind
(c) Observe that the linear bijection $\Psi: (\bigotimes_{i\in I_1} H_i) \otimes (\bigotimes_{j\in I_2} H_j) \to \bigotimes_{i\in I} H_i$ as in Remark \ref{rem:alt-con-ten-prod}(b) restricts to a surjection from $(\bigotimes_{i\in I_1}^\un H_i) \otimes (\bigotimes_{j\in I_2}^\un H_j)$ to $\bigotimes_{i\in I}^\un H_i$. 
Moreover, for any $u, u'\in \Pi_{i\in I_1} \sph(H_i)$ and $v,v'\in \Pi_{j\in I_2}  \sph(H_j)$, we have $(u,u')\sim (v,v')$ as elements in $\PI  \sph(H_i)$ if and only if $u\sim u'$ and $v\sim v'$. 
Thus, the argument is part (b) tells us that 
$$\left< (\otimes_{i\in I_1} u_i) \otimes (\otimes_{j\in I_2} v_j), 
(\otimes_{i\in I_1} u_i') \otimes (\otimes_{j\in I_2} v_j')\right>_{\phi_1}
\ =\  \langle \otimes_{i\in I_1} u_i, 
\otimes_{i\in I_1} u_i'\rangle_{\phi_1}
\langle \otimes_{j\in I_2} v_j, 
\otimes_{j\in I_2} v_j'\rangle_{\phi_1}$$
This shows that $\Psi\big|_{(\bigotimes_{i\in I_1}^\un H_i) \otimes (\bigotimes_{j\in I_2}^\un H_j)}$ is inner-product preserving. 
\end{prf}

\medskip

\emph{We denote by $\bar\bigotimes_{i\in I}^\mu H_i$ and $\bar\bigotimes_{i\in I}^{\phi_1} H_i$ the completions of $\bigotimes_{i\in I}^\mu H_i$ and ${\bigotimes}_{i\in I}^\un H_i$, respectively, under the norms induced by $\langle \cdot,\cdot \rangle_{\phi_1}$.} 

\medskip

\begin{eg}\label{eg:bad-expan}
If $H_i = \BC$ ($i\in I$), then the sesquilinear form $\langle \cdot, \cdot \rangle_{\phi_1}$ is not positive on the whole space $\BOI H_i$ since 
$\big\langle (\tsi 1/2 - \tsi 2), (\tsi 1/2 - \tsi 2)\big\rangle_{\phi_1} = -2$. 
\end{eg}

\medskip

Set $\PI^{\rm eu} H_i:= \{x\in \PI H_i: x_i\in \sph(H_i) \text{ except for a finite number of }i\}$ and 
$K$ be an inner-product space. 
A multilinear map $\Phi: \PI^{\rm eu} H_i \to K$ (i.e. $\Phi$ is coordinatewise linear) is said to be \emph{componentwise inner-product preserving} if for any $\mu,\nu\in \Omega^\un_{I;H}$,  
$$\left<\Phi(x), \Phi(y)\right> = \delta_{\mu,\nu} \PI\ \! \langle x_i, y_i \rangle 
\qquad (x\in \PI^\mu  H_i; y\in \PI^\nu  H_i)$$ 
where $\delta_{\mu,\nu}$ is the Kronecker delta. 

\medskip

\begin{thm}\label{thm:univ-prop-Hil-sp-ten-prod}
(a) $\bar\bigotimes_{i\in I}^{\phi_1} H_i \cong \bar\bigoplus_{\mu\in \Omega_{I;H}^\un}^{\ell^2} \bar\bigotimes_{i\in I}^{\mu} H_i$ canonically as Hilbert spaces. 

\smnoind
(b) $\Theta_H|_{\PI^{\rm eu}  H_i}: \PI^{\rm eu}  H_i \to \BOI^\un H_i$ is a componentwise inner-product preserving multilinear map. For any inner-product space $K$ and any componentwise inner-product preserving multilinear map $\Phi: \PI^{\rm eu}  H_i \to K$, there is a unique isometry $\ti \Phi: \BOI^\un H_i \to K$ such that $\Phi = \ti \Phi\circ \Theta_H|_{\PI^{\rm eu}  H_i}$.
\end{thm}
\begin{prf}
(a) Clearly, $\BOI^\un H_i = \sum_{\mu\in \Omega^\un_{I;H}} \BOI^{\mu} H_i$. 
Moreover, as in the proof of Lemma \ref{lem:phi0-innpro}(b), the two subspaces $\BOI^{\mu} H_i$ and $\BOI^{\nu} H_i$ are orthogonal if $\mu$ and $\nu$ are distinct elements in $\Omega^\un_{I;H}$. 
The rest of the argument is standard. 

\smnoind
(b) It is easy to see that $\Theta_H|_{\PI^{\rm eu}  H_i}$ is componentwise inner-product preserving. 
The uniqueness of $\ti \Phi$ follows from the fact that $\Theta_H(\PI^{\rm eu}  H_i)$ generates $\BOI^\un H_i$. 
To show the existence of  $\ti \Phi$, we first define a multilinear map $\Phi_0: \PI  H_i \to K$ by setting $\Phi_0 = \Phi$ on $\PI^{\rm eu}  H_i$ and 
$\Phi_0 = 0$ on $\PI  H_i \setminus \PI^{\rm eu}  H_i$. 
Let $\ti \Phi_0: \BOI H_i \to K$ be the induced linear map and set $\ti \Phi := \ti \Phi_0|_{\BOI^\un H_i}$. 
Suppose that $u,v\in \PI  \sph(H_i)$, $\xi\in \BOI^{u} H_i$ and $\eta\in \BOI^v H_i$. 
If $u\nsim v$, then $\langle \xi, \eta\rangle_{\phi_1} = 0 = 
\langle \ti \Phi(\xi), \ti \Phi(\eta)\rangle$. 
Otherwise, there exist $F\in \mathfrak{F}$ and $\xi_0, \eta_0\in \bigotimes_{i\in F} H_i$ such that $\xi = J^{u}_F(\xi_0)$, $\eta = J^{v}_F(\eta_0)$ and $u_i = v_i$ if $i\in I\setminus F$. 
In this case, 
$\langle \ti \Phi(\xi), \ti \Phi(\eta)\rangle
 = \langle \xi_0, \eta_0\rangle
 = \langle \xi, \eta\rangle_{\phi_1}$. 
\end{prf}

\medskip

\begin{eg}\label{eg:ten-Hil-sp}
Suppose that $\Phi$ and $\varphi$ are as in Corollary \ref{cor:boiut-C}, and $\{\delta_\mu\}_{\mu\in \Omega^\un_{I;\BC}}$ is the canonical orthonormal basis for $\ell^2\big(\Omega^\un_{I;\BC}\big)$. 
Note that $\Omega_{I;\BC}^\ut = \Omega_{I;\BC}^\un$ and consider the linear bijection $J: \BC[\Omega^\ut_{I;\BC}] \to \BC[\Omega^\un_{I;\BC}]$ given by $J(\lambda([\alpha]_\sim)) := \delta_{[\alpha]_\sim}$ ($\alpha\in \PI \BT$). 
By Example \ref{eg:u-ten-prod}(a) and Theorem \ref{thm:univ-prop-Hil-sp-ten-prod}(a), the map $J\circ \Phi$ induces a Hilbert space isomorphism $\hat\Phi:\bar\bigotimes_{i\in I}^{\phi_1} \BC \to \ell^2\big(\Omega^\un_{I;\BC}\big)$ such that $\hat\Phi(\tsi \beta_i) = \varphi(\beta)\delta_{[\beta]_\sim}$ ($\beta\in \PI\BT$). 
\end{eg}

\medskip

We would like to compare $\bar\bigotimes_{i\in I}^{\phi_1} H_i$ with the infinite directed product as defined in \cite{vN}, when $\{H_i\}_\ii$ is a family of Hilbert spaces. 
Let us first recall from \cite[Definition 3.3.1]{vN} that $x\in \PI   H_i$ is a \emph{$C_0$-sequence} if $\sum_{i\in I} \big| \|x_i\| - 1 \big|$ converges. 
As in \cite[Definition 3.3.2]{vN}, if $x$ and $y$ are $C_0$-sequences such that $\sum_{i\in I} \big| \langle x_i, y_i\rangle - 1 \big|$ converges, then we write $x \approx y$. 
Denote by $[x]_\approx$ the equivalence class of $x$ under $\approx$, and by  $\Gamma_{I;H}$ the set of all such equivalence classes (see \cite[Definition 3.3.3]{vN}). 

\medskip

Let $\prod\otimes_{i\in I} H_i$ be the infinite direct product Hilbert space as defined in \cite{vN}, and 
$\prod\otimes_\ii \ \!x_i$ be the element in $\prod{\otimes}_{i\in I} H_i$ corresponding to a $C_0$-sequence $x$ as in \cite[Theorem IV]{vN}.  
Notice that if $x\in \PI^{\rm eu}  H_i$, then $x$ is a $C_0$-sequence, and we have a multilinear map 
$$\Upsilon: \PI^{\rm eu}  H_i\ \longrightarrow\  \prod{\otimes}_{i\in I} H_i.$$ 
On the other hand, for any $\mathfrak{C}\in \Gamma_{I;H}$, we denote $\prod\otimes_{i\in I}^\mathfrak{C} H_i$ to be the closed subspace of $\prod\otimes_{i\in I} H_i$ generated by $\{ \prod\otimes_\ii\ \! x_i : x\in \mathfrak{C}\}$
(see \cite[Definition 4.1.1]{vN}). 

\medskip

\begin{prop}\label{prop:cp-vN}
Let $\{H_i\}_{i\in I}$ be a family of Hilbert spaces. 

\smnoind
(a) $[x]_\sim \mapsto [x]_\approx$ ($x\in \PI  \sph(H_i)$) gives a well-defined surjection $\kappa_H: \Omega^\un_{I;H} \to \Gamma_{I;H}$. 
Moreover, for any $x,y\in \PI  \sph(H_i)$, there is a bijection between $\kappa_H^{-1}([x]_\approx)$ and $\kappa_H^{-1}([y]_\approx)$.

\smnoind
(b) There exists a linear map $\ti \Upsilon: \BOI^\un H_i \to \prod{\otimes}_{i\in I} H_i$ such that $\Upsilon = \ti \Upsilon \circ \Theta_H |_{\PI^{\rm eu} H_i}$ and $\ti \Upsilon\mid_{\BOI^\mu H_i}$ extends to a Hilbert space isomorphism $\ti\Upsilon^\mu: \bar\bigotimes_{i\in I}^{\mu} H_i \to \prod\otimes_{i\in I}^{\kappa_H(\mu)} H_i$ ($\mu\in \Omega^\un_{I;H}$). 
\end{prop}
\begin{prf}
(a) Clearly, if $x\sim z$, then $x\approx z$ and $\kappa_H$ is well-defined. 
\cite[Lemma 3.3.7]{vN} tells us that $\kappa_H$ is surjective. 
Furthermore, there exists a unitary $u_i\in \CL(H_i)$ such that $u_ix_i = y_i$ ($i\in I$), and $[u_i]_{i\in I}$ induces the required bijective correspondence in the second statement. 

\smnoind
(b) By the argument of Theorem \ref{thm:univ-prop-Hil-sp-ten-prod}(b), one can construct a linear map $\ti \Upsilon$ such that $\Upsilon = \ti \Upsilon \circ \Theta_H |_{\PI^{\rm eu}  H_i}$. 
By the argument of part (a), we see that $\ti \Upsilon\left(\BOI^{[u]_\sim} H_i\right) \subseteq \prod\otimes_{i\in I}^{[u]_\approx} H_i$ ($u\in \PI  \sph(H_i)$). 
Furthermore, by Lemma \ref{lem:phi0-innpro}(a), Proposition \ref{prop:ten-inn-prod}(b) and \cite[Theorem IV]{vN}, we see that $\ti \Upsilon|_{\BOI^{[u]_\sim} H_i}$ is an isometry. 
Finally, $\ti \Upsilon|_{\BOI^{[u]_\sim} H_i}$ has dense range (by \cite[Lemma 4.1.2]{vN}). 
\end{prf}

\medskip

Notice that $\ti \Upsilon$ is, in general, unbounded but 
Remark \ref{rem:cp-inf-direct-prod}(b) below tells us that $\bar\bigotimes_\ii^{\phi_1} H_i$ is a ``natural dilation'' of $\prod \tsi H_i$. 
On the other hand, Remark \ref{rem:cp-inf-direct-prod}(d) says that it is possible to construct $\prod \tsi H_i$ in a similar way as that for $\bar\bigotimes_\ii^{\phi_1} H_i$. 
Note however, that the construction of $\bar\bigotimes_\ii^{\phi_1} H_i$ is totally algebraical and $\bar\bigotimes_\ii^{\phi_1} H_i$ itself seems to be more natural (see Theorem \ref{thm:ten-prod-hil-cT-mod} and Example \ref{eg:ten-prod-rep-BC} below). 

\medskip

\begin{rem}\label{rem:cp-inf-direct-prod}
Suppose that $\{H_i\}_{i\in I}$ is a family of Hilbert spaces. 

\smnoind
(a) $\sim$ and $\approx$ are different even in the case when $I = \mathbb{N}$ and $H_i = \BC$ ($i\in \BN$) because one can find $x,y\in \Pi_{i\in \BN} \BT$ with $x_i \neq y_i$ for all $i\in \BN$ but $\sum_{i=1}^\infty \big|\langle x_i, y_i \rangle - 1\big|$ converges. 
In fact, $\kappa_H^{-1}([x]_\approx)$ is an infinite set. 

\smnoind
(b) By \cite[Lemma 4.1.1]{vN}, we have 
\begin{equation*}\label{eqt:decomp-inf-dir-prod}
\prod{\otimes}_{i\in I} H_i\ =\ {\bar\bigoplus}_{\mathfrak{C}\in \Gamma_{I;H}}^{\ell^2} \prod{\otimes}_{i\in I}^\mathfrak{C} H_i.
\end{equation*}
Therefore, Theorem \ref{thm:univ-prop-Hil-sp-ten-prod}(a) and  Proposition \ref{prop:cp-vN} tell us that for a fixed $\gamma_0\in \Gamma_{I;H}$, one has a canonical Hilbert space isomorphism $${\bar\bigotimes}_{i\in I}^{\phi_1} H_i
\ \cong\ \ell^2\big(\kappa_H^{-1}(\gamma_0)\big) \bar \otimes \big(\prod {\otimes}_{i\in I} H_i\big).$$ 

\smnoind
(c) For each $i\in I$, let $K_i$ be an inner-product space such that $H_i$ is the completion of $K_i$. 
Then $\bar\bigotimes_{i\in I}^{\phi_1} K_i$ is, in general, not canonically isomorphic to $\bar\bigotimes_{i\in I}^{\phi_1} H_i$ because $\Omega_{I;K}^\un \subsetneq \Omega_{I;H}^\un$ if $K_i\subsetneq H_i$ for an infinite number of $i\in I$. 
On the other hand, if $I$ is countable, for any $x\in \PI  \sph(H_i)$, there exists $y\in \PI  \sph(K_i)$ such that $x\approx y$. 
This shows that the restriction, $\kappa_{H;K}$, of $\kappa_H$ on $\Omega^\un_{I;K}$ is also a surjection onto $\Gamma_{I;H}$.  
However, we do not know if the cardinality of $\kappa_{H;K}^{-1}(\mathfrak{C})$ are the same for different $\mathfrak{C}\in \Gamma_{I;H}$. 

\smnoind
(d) If $\phi_0$ is as in Example \ref{eg:multi-linear}(b), it is easy to see that 
$$\langle \prod \otimes u_i, \prod \otimes v_i\rangle
 = \phi_0\big(\langle \tsi u_i, \tsi v_i\rangle_{\BC^{\otimes I}} \big)
\qquad (u,v\in \PI^\un H_i).$$
Thus, the sesquilinear form $\phi_0\big(\langle\cdot, \cdot\rangle_{\BC^{\otimes I}} \big)$ produces $\prod \otimes H_i$. 
If one wants a self-contained alternative construction for $\prod \otimes H_i$, one needs to establish the positivity of $\phi_0\big(\langle\cdot, \cdot\rangle_{\BC^{\otimes I}} \big)$, which can be reduced to showing the positivity when all $H_i$ are of the same finite dimension. 
\end{rem}

\medskip

In the remainder of this section, we show that $\BOI^\un H_i$ can be completed into a $C^*(\Omega^\ut_{I;\BC})$-module, which gives many pre-inner products on $\BOI^\un H_i$ including $\langle\cdot , \cdot\rangle_{\phi_1}$. 
In the following, we use the convention that the $A$-valued inner-product of an inner-product $A$-module is $A$-linear in the first variable (where $A$ is a pre-$C^*$-algebra). 
On the other hand, we recall that if $G$ is a group and $\lambda_g$ is the canonical image of $g$ in $\BC[G]$, the map ${\sum}_{g\in G} \alpha_g \lambda_g \mapsto \alpha_e$ ($\alpha_g\in \BC$), where $e\in G$ is the identity, extends to a faithful tracial state $\chi_G$ on $C^*(G)$. 

\medskip 

\begin{thm}\label{thm:ten-prod-hil-cT-mod}
(a) There exists an inner-product $\BC[\Omega^\ut_{I;\BC}]$-module structure on $\BOI^\un H_i$. 
If $\bar\bigotimes_{i\in I}^{\rm mod} H_i$ is the Hilbert $C^*(\Omega^\ut_{I;\BC})$-module given by the completion of this $\BC[\Omega^\ut_{I;\BC}]$-module, we have a canonical Hilbert space isomorphism
\begin{equation}\label{eqt:ten-phi=mod-ten}
{\bar\bigotimes}_{i\in I}^{\phi_1} H_i\ \cong\ \big({\bar\bigotimes}_{i\in I}^{\rm mod} H_i\big) \bar\otimes_{\chi_{\Omega^\ut_{I;\BC}}} \BC.
\end{equation}

\smnoind
(b) If $G\subseteq \Omega^\ut_{I;\BC}$ is a subgroup and $\mathcal{E}_G: C^*(\Omega^\ut_{I;\BC}) \to C^*(G)$ is the canonical conditional expectation, there is an inner-product $\BC[G]$-module structure on $\BOI^\un H_i$, whose completion coincides with the Hilbert $C^*(G)$-module $\big({\bar\bigotimes}_{i\in I}^{\rm mod} H_i\big) \bar\otimes_{\mathcal{E}_G} C^*(G)$.
\end{thm}
\begin{prf}
(a) Clearly, $\BOI^\un H_i$ is a $\BC^{\otimes I}_\ut$-submodule of the $\BC^{\otimes I}$-module $\BOI H_i$ (see Proposition \ref{prop:exist-ten-prod}(c)).
Moreover, 
one has a linear ``truncation'' 
$E$ from $\BC^{\otimes I} = \big({\bigoplus}_{\omega\in \Omega_{I;\BC}\setminus \Omega^\ut_{I;\BC}} \BOI^\omega \BC\big) \oplus \BC^{\otimes I}_\ut$ to $\BC^{\otimes I}_\ut$ sending $(\alpha, \beta)$ to $\beta$. 
Define
$$\langle \xi, \eta \rangle_{\BC^{\otimes I}_\ut}
 := E\big( \langle \xi, \eta \rangle_{\BC^{\otimes I}} \big)
\qquad \big(\xi, \eta\in \BOI^\un H_i\big),$$
which is a Hermitian $\BC^{\otimes I}_\ut$-form because by \eqref{eqt:grading}, we have 
$$E(ab) = E(a)b \quad \text{and}\quad  E(a^*) = E(a)^* 
\qquad (a\in \BC^{\otimes I};b\in \BC^{\otimes I}_\ut).$$ 
For any $u,v\in \PI  \sph(H_i)$, we write $u\sim_\s v$ if there exists $\beta\in \PI \BT$ such that $u_i = \beta_i v_i$ e.f. 
Then $\sim_\s$ is an equivalence relation on $\PI  \sph(H_i)$ satisfying
\begin{equation}\label{eqt:sims-equiv}
u\sim_\s v\quad \text{if and only if}\quad \langle \tsi u_i, \tsi v_i \rangle_{\BC^{\otimes I}} \in \BC^{\otimes I}_\ut.
\end{equation}

Let $\Phi$ and $\varphi$ be as in Corollary \ref{cor:boiut-C}.
Suppose that $\xi = \sum_{k=1}^n \alpha_k \tsi u^{(k)}_i$ with $\alpha_1,...,\alpha_n\in \BC$ and $u^{(1)}, ..., u^{(n)}\in \PI  \sph(H_i)$. 
We first show that $\Phi(\langle \xi, \xi \rangle_{\BC^{\otimes I}_\ut}) \in C^*(\Omega^\ut_{I;\BC})_+$.
As in the proof of Lemma \ref{lem:phi0-innpro}(b), it suffices to consider the case when $u^{(k)} \sim_\s u^{(1)}$ for any $k \in \{1,...,n\}$ (because of Relation \eqref{eqt:sims-equiv}). 
Let $F\in \mathfrak{F}$ and $\beta^{(1)},...,\beta^{(n)}\in \PI \BT$ such that $u^{(k)}_i = \beta^{(k)}_iu^{(1)}_i$ ($i\in I\setminus F$; $k=1,...,n$). 
For any $k,l\in \{1,...,n\}$, we have
$$
\Phi\big((\Pi_{i\in F} \langle u^{(k)}_i, u^{(l)}_i \rangle_i) (\otimes_{i\in I\setminus F}\ \!\beta_i^{(k)}
\overline{\beta_i^{(l)}})\big) 
\ =\ \langle \ti\varphi_F(u^{(k)}), \ti\varphi_F(u^{(l)}) \rangle_F,$$ 
where $\ti\varphi_F(u^{(k)}) := \big(\varphi(\beta^{(k)}) \Pi_{i\in F} \beta^{(k)}_i\big)^{-1} (\otimes_{i\in F}\ \! u^{(k)}_i)\otimes \lambda_{[\beta^{(k)}]_\sim}$ and 
$\langle \cdot, \cdot \rangle_F$ is the canonical $\BC[\Omega^\ut_{I;\BC}]$-valued inner-product on $(\bigotimes_{i\in F} H_i) \otimes \BC[\Omega^\ut_{I;\BC}]$. 
Therefore, 
\begin{equation*}\label{eqt:mod-pos}
\Phi(\langle \xi ,\xi \rangle_{\BC^{\otimes I}_\ut})
\ = \ \left\langle {\sum}_{k=1}^n \alpha_k \ti\varphi_F(u^{(k)}), {\sum}_{k=1}^n \alpha_k \ti\varphi_F(u^{(k)}) \right\rangle_F
\ \geq \ 0.
\end{equation*}

Next, we show that $\chi_{\Omega^\ut_{I;\BC}}\circ \Phi\circ E = \phi_1$. 
Let $\alpha\in \PI \BC^\times$. 
If $\alpha\nsim 1$, then $\chi_{\Omega^\ut_{I;\BC}}\circ \Phi\circ E(\tsi \alpha_i) = 0$ (as $\Phi(E(\tsi \alpha_i))\notin \BC\cdot \lambda_{[1]_\sim}\setminus \{0\}$, whether or not $[\alpha]_\sim\in \Omega^\ut_{I;\BC}$) and we also have $\phi_1(\tsi \alpha_i) = 0$.
If $\alpha \sim 1$, then $\tsi \alpha_i = (\PI \alpha_i)(\tsi 1) = (\PI \alpha_i)\lambda_{[1]_\sim}$, which implies that 
$\chi_{\Omega^\ut_{I;\BC}}(\Phi(\tsi \alpha_i)) 
= \PI \alpha_i 
= \phi_1(\tsi \alpha_i)$.

Thus, we have 
\begin{equation}\label{eqt:cp-chi-phi0}
\chi_{\Omega^\ut_{I;\BC}}\big(\Phi(\langle \xi, \eta \rangle_{\BC^{\otimes I}_\ut})\big)
 = \langle \xi, \eta \rangle_{\phi_1}
\qquad \big(\xi, \eta\in \BOI^\un H_i\big).
\end{equation}
As a consequence, if $\Phi(\langle \xi, \xi \rangle_{\BC^{\otimes I}_\ut}) = 0$, we know from Lemma \ref{lem:phi0-innpro}(b) that $\xi = 0$. 
This gives an inner-product $\BC[\Omega^\ut_{I;\BC}]$-module structure on $\BOI^\un H_i$. 
Furthermore, the required isomorphism $\bar\bigotimes_{i\in I}^{\phi_1} H_i \cong (\bar\bigotimes_{i\in I}^{\rm mod} H_i) \bar \otimes_{\chi_{\Omega^\ut_{I;\BC}}} \BC$ also follows from 
\eqref{eqt:cp-chi-phi0}. 

\smnoind
(b) Since $\BOI^\un H_i$ is a $\BC[G]$-module (under the identification of $\BC[G]$ with $\bigoplus_{\omega\in G} \BOI^\omega \BC$ under the $^*$-isomorphism $\Phi$ in Corollary \ref{cor:boiut-C}), every element in $(\BOI^\un H_i)\otimes_{\BC[G]} \BC[G]$ is of the form $\xi\otimes_{\BC[G]} 1$ for some $\xi\in \BOI^\un H_i$. 
Moreover, if $\xi,\eta\in \BOI^\un H$, then 
\begin{equation}\label{eqt:G-mod}
\langle \xi\otimes_{\BC[G]} 1, \eta\otimes_{\BC[G]} 1\rangle_{({\bar\bigotimes}_{i\in I}^{\rm mod} \BC)\bar\otimes_{\mathcal{E}_G} C^*(G)} 
\ =\ \mathcal{E}_G(\Phi(\langle \xi, \eta\rangle_{\BC^{\otimes I}_\ut})) 
\ =\ \Phi(E_G(\langle \xi, \eta\rangle_{\BC^{\otimes I}})), 
\end{equation}
where $E_G$ is the linear ``truncation'' map from $\BC^{\otimes I}$ to $\bigoplus_{\omega\in G} \BOI^\omega \BC$ defined as in part (a). 
Therefore, $\Phi(E_G(\langle \cdot, \cdot\rangle_{\BC^{\otimes I}}))$ is a positive Hermitian $\BC[G]$-form on $\BOI^\un H_i$. 
Obviously, $\chi_{\Omega^\ut_{I;\BC}} = \chi_G\circ \mathcal{E}_G$, and by \eqref{eqt:cp-chi-phi0}, 
$$\chi_G(\Phi(E_G(\langle \xi, \eta\rangle_{\BC^{\otimes I}})))
 = \chi_{\Omega^\ut_{I;\BC}}(\Phi(\langle \xi, \eta\rangle_{\BC^{\otimes I}_\ut}))
 = \langle \xi, \eta\rangle_{\phi_1}
\qquad \big(\xi,\eta\in \BOI^\un H\big).$$ 
This implies that $\Phi(E_G(\langle \cdot, \cdot\rangle_{\BC^{\otimes I}}))$ is non-degenerate (since $\langle\cdot, \cdot\rangle_{\phi_1}$ is non-degenerate by Lemma \ref{lem:phi0-innpro}(b)). 
Now, Equation \eqref{eqt:G-mod} tells us that the Hilbert $C^*(G)$-module $\big({\bar\bigotimes}_{i\in I}^{\rm mod} H_i\big)\bar\otimes_{\mathcal{E}_G} C^*(G)$ is the completion of $\BOI^\un H_i$ under the norm induced by the $\BC[G]$-valued inner-product $\Phi(E_G(\langle \cdot, \cdot\rangle_{\BC^{\otimes I}}))$. 
\end{prf}

\medskip

Let $\{e\}$ be the trivial subgroup of $\Omega^\ut_{I;\BC}$. 
Since one can identify $E_{\{e\}}$ with $\phi_1$ (through the argument of Theorem \ref{thm:ten-prod-hil-cT-mod}(b)), one has 
$${\bar\bigotimes}_{i\in I}^{\phi_1} H_i \ \cong \ \big({\bar\bigotimes}_{i\in I}^{\rm mod} H_i\big)\bar\otimes_{\mathcal{E}_{\{e\}}} \BC.$$

\medskip

\begin{rem}\label{rem:hil-mod}
(a) For any subgroup $G\subseteq \Omega^\ut_{I;\BC}$ and any faithful state $\varphi$ on $C^*(G)$, 
the Hilbert space 
$$\Big(\big({\bar\bigotimes}_\ii^{\rm mod} H_i\big)\bar\otimes_{\mathcal{E}_G} C^*(G)\Big) \bar\otimes_\varphi \BC$$ 
induces an inner-product on $\BOI^\un H_i$. 

\smnoind
(b) If $x\in \PI^0 \BC$ (see Example \ref{eg:multi-linear}(b)), then $\sup_{i\in I} \abs{x_i} < \infty$.
This, together with the surjectivity of $\kappa_\BC$ (see Proposition \ref{prop:cp-vN}(a)), tells us that $\Gamma_{I;\BC}$ is a group under the multiplication: $[x]_\approx\cdot [y]_\approx := [xy]_\approx$ (where $(xy)_i := x_iy_i$ for any $i\in I$). 
Moreover, $\kappa_\BC: \Omega^\ut_{I;\BC} = \Omega^\un_{I;\BC} \to \Gamma_{I;\BC}$ is a group homomorphism, which induces a surjective $^*$-homomorphism $\bar \kappa_\BC: C^*(\Omega^\ut_{I;\BC}) \to C^*(\Gamma_{I;\BC})$.  

\smnoind
(c) It is natural to ask whether $\big((\bar\bigotimes_\ii^{\rm mod} H_i) \bar\otimes_{\bar\kappa_\BC} C^*(\Gamma_{I;\BC})\big)\bar\otimes_{\chi_{\Gamma_{I;\BC}}}\BC$ is isomorphic to $\prod \tsi H_i$ canonically. 
Unfortunately, it is not the case. 
In fact, for any $x,y\in \PI^\un H_i$, we denote $x \approx_\BT y$ if there exists $\alpha\in \PI \BT$ with $\alpha\approx 1$ such that $x_i = \alpha_i y_i$ e.f.
It is easy to check that $\approx_\BT$ is an equivalent relation standing strictly between $\sim$ and $\approx$ in general. 
Moreover, one has 
$$\big< ((\tsi x_i)\otimes_{\bar\kappa_\BC} 1)\otimes_{\chi_{\Gamma_{I;\BC}}} 1, ((\tsi y_i)\otimes_{\bar\kappa_\BC} 1)\otimes_{\chi_{\Gamma_{I;\BC}}} 1\big>
 = 0
\quad \text{whenever } x\not\approx_\BT y,$$ 
while $\big< \prod \tsi x_i, \prod \tsi y_i\big> = 0$ whenever $x\not\approx y$. 
Note however, that if all $H_i = \BC$, then $\approx_\BT$ and $\approx$ coincide, and one can show that the two Hilbert spaces  
$\big((\bar\bigotimes_\ii^{\rm mod} \BC) \bar\otimes_{\kappa_\BC} C^*(\Gamma_{I;\BC})\big)\bar\otimes_{\chi_{\Gamma_{I;\BC}}}\BC$ and $\prod \tsi \BC$ coincide canonically. 
\end{rem}

\medskip

\begin{eg}\label{eg:ten-Hil-mod}
(a) It is clear that $\bar\bigotimes_\ii^{\rm mod} \BC = C^*(\Omega^\ut_{I;\BC})$. 
For any state $\varphi$ on $C^*(\Omega^\ut_{I;\BC})$, the Hilbert space $(\bar\bigotimes_\ii^{\rm mod} \BC)\bar\otimes_\varphi \BC$ is the GNS construction of $\varphi$. 

\smnoind
(b) If $G$ is a subgroup of $\Omega^\ut_{I;\BC}$, we have
$$\big({\bar\bigotimes}_{i\in I}^{\rm mod} \BC\big)\bar\otimes_{\mathcal{E}_G} C^*(G)
\ \cong\ \ell^2(\Omega^\ut_{I;\BC}/G) \bar\otimes C^*(G).$$
In fact, let $q: \Omega^\ut_{I;\BC}\to \Omega^\ut_{I;\BC}/G$ be the quotient map and $\sigma: \Omega^\ut_{I;\BC}/G \to \Omega^\ut_{I;\BC}$ be a cross-section. 
One has a bijection from $\Omega^\ut_{I;\BC}$ to $(\Omega^\ut_{I;\BC}/G)\times G$ sending $\omega$ to $(q(\omega),  \sigma(q(\omega)^{-1})\omega)$. 
This induces a bijective linear map $\Delta : \BC[\Omega^\ut_{I;\BC}] \to \bigoplus_{\Omega^\ut_{I;\BC}/G} \BC[G]$ such that 
for any $\omega\in \Omega^\ut_{I;\BC}$ and $\varepsilon \in \Omega^\ut_{I;\BC}/G$, 
$$\Delta(\lambda_{\omega})_\varepsilon
\ := \ \begin{cases}
\lambda_{\sigma(\varepsilon^{-1})\omega} \ \ & \text{if }q(\omega)= \varepsilon\\
0 & \text{otherwise}.  
\end{cases}$$
Let $\Phi: \BOI^\un \BC = \BC^{\otimes I}_\ut \to \BC[\Omega^\ut_{I;\BC}]$ and $\varphi:\PI \BT\to \BT$ be as in Corollary \ref{cor:boiut-C}. 
Suppose that $\alpha,\beta\in \PI \BC^\times$. 
If $[\alpha\beta^{-1}]_\sim$ does not belong to $G$, then $E_G(\langle \tsi \alpha_i, \tsi \beta_i\rangle_{\BC^{\otimes I}}) = 0$, and 
$$\big\langle \Delta\circ\Phi\big(\tsi \alpha_i\big), \Delta\circ\Phi(\tsi \beta_i\big) \big\rangle_{\bigoplus_{\Omega^\ut_{I;\BC}/G}^{\ell^2} \BC[G]}
\ = \ 0.$$
On the other hand, if $[\alpha\beta^{-1}]_\sim\in G$, then 
\begin{eqnarray*}
\big\langle \Delta\circ\Phi\big(\tsi \alpha_i\big), \Delta\circ\Phi(\tsi \beta_i\big) \big\rangle_{\bigoplus_{\Omega^\ut_{I;\BC}/G}^{\ell^2} \BC[G]}
& = & \varphi(\alpha \beta^{-1}) \lambda_{[\alpha\beta^{-1}]_\sim}\\ 
& = & \Phi(\tsi \alpha_i\beta^{-1}_i)
\  = \ \Phi(E_G(\langle \tsi \alpha_i, \tsi \beta_i\rangle_{\BC^{\otimes I}})).
\end{eqnarray*}
This shows that $\Delta\circ \Phi$ is an inner-product $\BC[G]$-module isomorphism from $\BOI^\un \BC$ (equipped with the inner-product $\BC[G]$-module structure as in Theorem  \ref{thm:ten-prod-hil-cT-mod}(b)) onto $\bigoplus_{\Omega^\ut_{I;\BC}/G}^{\ell^2} \BC[G]$. 
\end{eg}

\medskip

\section{Tensor products of $^*$-representations of $^*$-algebras}

\medskip

\emph{In this section, $\{(A_i, H_i, \Psi_i)\}_{i\in I}$ is a family of unital $^*$-representations, in the sense that $A_i$ is a unital $^*$-algebra, $H_i$ is a Hilbert space and $\Psi_i: A_i \to \CL(H_i)$ is a unital $^*$-homomorphism ($i\in I$).} 

\medskip

Suppose that $\Psi_0:=\TBOIM\Psi_i: \BOI A_i \to L(\BOI H_i)$ is the map as in Proposition \ref{prop:exist-ten-prod}(c). 
It is easy to check that 
\begin{equation}\label{eqt:ten-psi-inv}
\big< \Psi_0 (a) \xi, \eta\big>_{\BC^{\otimes I}}
=\big< \xi, \Psi_0(a^*) \eta\big>_{\BC^{\otimes I}} 
\qquad \big(a\in {\bigotimes}_{i\in I} A_i;\xi,\eta\in {\bigotimes}_{i\in I} H_i\big).
\end{equation}
Furthermore, one has the following result (which is more or less well-known). 

\medskip

\begin{prop}\label{prop:ten-prod-st-alg}
For any $\mu\in \Omega^\un_{I;H}$, the map $\TBOIM\Psi_i$ induces a unital $^*$-representation ${\bigotimes}_{i\in I}^\mu\Psi_i: \bigotimes_{i\in I}^e A_i \to \CL(\bar\bigotimes_{i\in I}^\mu H_i)$. 
If, in addition, all $\Psi_i$ are injective, then so is $\BOI^\mu\Psi_i$. 
\end{prop}

\medskip

Consequently, one has a unital $^*$-representation of $\BOI^e A_i$ on the Hilbert space $\bar\bigotimes^{\phi_1}_{i\in I} H_i$. 
However, it seems impossible to extend it to a unital $^*$-representation of $\BOI A_i$ on $\bar\bigotimes^{\phi_1}_{i\in I} H_i$.
The biggest $^*$-subalgebra $\BOI A_i$ that we can think of, for which such extension is possible, is the subalgebra $\BOI^\ut A_i$. 
Example \ref{eg:ten-prod-rep-BC}(a) also tells us that it is probably the right subalgebra to consider. 

\medskip

Let us digress a little bit and give another $^*$-representation of $\BOI^\ut A_i$, which is a direct consequence of Proposition \ref{prop:ten-prod-st-alg}, Theorem \ref{thm:boiut-tw-cr-pd}(a) and \cite[Theorem 4.1]{BS} (it is not hard to verify that the representation as given in \cite[Theorem 4.1]{BS} is injective when $\BOI^{\mu}\Psi_i$ is injective). 
Note however, that such a $^*$-representation is not canonical since it depends on the choices a cross-section $c:\Omega^\ut_{I;A} \to \PI  U_{A_i}$ (see Remark \ref{rem:ut-prod}(a)). 

\medskip

\begin{cor}
Suppose that $\Psi_i$ are injective.
For any $\mu\in \Omega^\un_{I;H}$, the injection $\BOI^{\mu}\Psi_i$ induces an injective unital $^*$-representation of 
$\BOI^\ut A_i$ on $(\bar\bigotimes_{i\in I}^\mu H_i)\otimes \ell^2(\Omega^\ut_{I;A})$.
\end{cor}

\medskip 

Let us now go back to the discussion of the tensor product type representation of $\BOI^\ut A_i$. 
Observe that $\{\Psi_i\}_\ii$ induces a canonical action $\alpha^\Psi: \Omega^\ut_{I;A}\times \Omega^\un_{I;H} \to \Omega^\un_{I;H}$. 
For simplicity, we will denote $\alpha^\Psi_\omega(\mu)$ by $\omega\cdot \mu$ ($\omega\in \Omega^\ut_{I;A}; \mu\in \Omega^\un_{I;H}$). 

\smnoind

\medskip

\begin{thm}\label{thm:inf-ten-c-st-alg}
(a) The map $\TBOIM \Psi_i$ induces a unital $^*$-representation $\bigotimes_{i\in I}^{\phi_1} \Psi_i: \BOI^\ut A_i \to \CL\big( \bar\bigotimes_{i\in I}^{\phi_1} H_i \big)$. 

\smnoind
(b) $\big( \bar\bigotimes_{i\in I}^{\phi_1} H_i, ({\bigotimes}_{i\in I}^{\phi_1} \Psi_i)|_{\BOI^eA_i}\big) = \bigoplus_{\mu\in \Omega_{I;H}^\un} \big(\bar\bigotimes_{i\in I}^{\mu} H_i, {\bigotimes}_{i\in I}^{\mu} \Psi_i \big)$.

\smnoind
(c) If all $\Psi_i$ are injective, then so is $\BOI^{\phi_1}\Psi_i$. 
\end{thm}
\begin{prf}
(a) Set $\Psi_0:=\TBOIM \Psi_i$. 
For any $\mu\in \Omega^\un_{I;H}$, $\omega\in \Omega^\ut_{I;A}$ and $a\in \PI^\omega  A_i$, it is clear that 
\begin{equation}\label{eqt:Psi-otimes-mu}
\Psi_0(\tsi a_i)\big(\BOI^\mu H_i \big)\ \subseteq\ \BOI^{\omega\cdot \mu} H_i. 
\end{equation}
Suppose that $u\in \omega$ and $F\in \mathfrak{F}$ such that $a_i = u_i$ for $i\in I\setminus F$.
If $\xi= J_{F'}^x(\xi_0)$ where $x\in \mu$, $F'\in \mathfrak{F}$ with $F\subseteq F'$ and $\xi_0\in \bigotimes_{i\in F'} H_i$, then 
$$\left< \Psi_0(\tsi a_i) \xi,  
\Psi_0(\tsi a_i) \xi \right>_{\BC^{\otimes I}} 
\ =\ \big< \big({\bigotimes}_{i\in F} \Psi_i(a_i)\otimes \id\big)\xi_0, \big({\bigotimes}_{i\in F} \Psi_i(a_i)\otimes \id\big)\xi_0 \big> (\tsi 1).$$
This means that $\Psi_0(\tsi a_i)$ is bounded on $\big(\BOI^\un H_i,\langle \cdot, \cdot\rangle_{\phi_1}\big)$ (see Theorem \ref{thm:univ-prop-Hil-sp-ten-prod}(a) and Proposition \ref{prop:ten-inn-prod}(b)) and produces a unital homomorphism $\BOI^{\phi_1} \Psi_i: \BOI^\ut A_i \to \CL\big( \bar\bigotimes_{i\in I}^{\phi_1} H_i \big)$.
Now, Relation \eqref{eqt:ten-psi-inv} tells us that $\BOI^{\phi_1} \Psi_i$ preserves the involution. 

\smnoind
(b) This part follows directly from the argument of part (a). 

\smnoind
(c) Set $\Psi := \BOI^{\phi_1} \Psi_i$. 
Suppose that $v^{(1)}, ..., v^{(n)}\in \PI  U_{A_i}$ are mutually inequivalent elements, $F\in \mathfrak{F}$, $b^{(1)},...,b^{(n)}\in \bigotimes_{i\in F} A_i$ and $a^{(k)} := J_F^{v^{(k)}}(b^{(k)})$ ($k=1,...,n$) such that  
\begin{equation*}\label{eqt:ker-Psi}
\Psi\big({\sum}_{k=1}^n a^{(k)}\big)\ =\ 0.
\end{equation*}
By induction, it suffices to show that $a^{(1)} = 0$. 

By replacing $a^{(k)}$ with $(v^{(1)})^{-1}a^{(k)}$ if necessary, we may assume that $v^{(1)}_i = e_i$ ($\ii$). 
If $n=1$, we take an arbitrary $\xi\in \PI \sph(H_i)$. 
If $n >1$, we claim that there exists $\xi\in \PI\sph(H_i)$ such that 
\begin{equation}\label{eqt:strong-faith}
\xi\ \nsim\ [V^{(k)}_i\xi_i]_\ii
\qquad (k=2,...,n),
\end{equation}
where $V^{(k)}_i := \Psi_i(v^{(k)}_i)$.
In fact, if $k\in \{2,...,n\}$ and  $i\in I^k := \{i\in I: v^{(k)}_i \neq e_i\}$ (which is an infinite set), the subset $\sph (H_i) \cap \ker (V^{(k)}_i - \id_{H_i})$ is nowhere dense in $\sph (H_i)$ as $\ker (V^{(k)}_i - \id_{H_i})$ is a proper closed subspace of $H_i$ (note that $\Psi_i$ is injective). 
For any $i\in I$, we consider $N_i := \{k\in\{2,...,n\}: i\in I^k\}$. 
By the Baire Category Theorem, for every $i\in I$, one can choose $\xi_i\in \sph(H_i) \setminus \bigcup_{k\in N_i} \ker (V^{(k)}_i - \id_{H_i})$. 
Now, $\xi:=[\xi_i]_\ii$ will satisfy Relation \eqref{eqt:strong-faith}. 

Since $\Psi(a^{(1)})\big(\BOI^\xi H_i\big) \subseteq \BOI^\xi H_i$ and 
$$\BOI^\xi H_i\ \cap\ {\sum}_{k=2}^n\Psi(a^{(k)})\big(\BOI^\xi H_i\big)\ =\ \{0\}$$ 
(because of Theorem \ref{thm:bas-inf-ten} as well as \eqref{eqt:Psi-otimes-mu} and \eqref{eqt:strong-faith}), 
we have $\Psi(a^{(1)})|_{\BOI^\xi H_i} = 0$. 
Therefore, part (b) and Proposition \ref{prop:ten-prod-st-alg} tells us that $a^{(1)} = 0$. 
\end{prf}

\medskip

\begin{rem}\label{rem-act-OA-OH}
(a) By the argument of Theorem \ref{thm:inf-ten-c-st-alg}(c), if all $\Psi_i$ are injective, then $\alpha^\Psi$ is \emph{strongly faithful} in the sense that for any finite subset $F\subseteq \Omega^\ut_{I;A} \setminus \{e\}$, there exists $\mu\in \Omega^\un_{I;H}$ with $\omega\cdot \mu \neq \mu$ ($\omega\in F$). 

\smnoind
(b) If $y, z\in \PI  H_i$ are $C_0$-sequences and $u,v \in \PI  U_{A_i}$, then 
\begin{equation}\label{eqt:defn-ti-alpha}
y \approx z \quad \text{if and only if} \quad [\Psi_i(u_i)y_i]_{i\in I}\approx [\Psi_i(u_i)z_i]_{i\in I}
\end{equation}
and $[\Psi_i(u_i)y_i]_{i\in I}\approx [\Psi_i(v_i)y_i]_{i\in I}$ whenever $u\sim v$. 
Thus, $\{\Psi_i\}_{i\in I}$ induces an action 
$\ti \alpha^\Psi: \Omega^\ut_{I;A}\times \Gamma_{I;H} \to \Gamma_{I;H}$. 
Again, we write $\omega\cdot \gamma$ for $\ti \alpha^\Psi_\omega(\gamma)$ ($\omega\in \Omega^\ut_{I;A}; \gamma\in\Gamma_{I;A}$).  
The map $\kappa_H$ in Proposition \ref{prop:cp-vN}(a) is \emph{equivariant} in the sense that  $\kappa_H\circ \alpha^\Psi_\omega = \ti\alpha^\Psi_\omega\circ \kappa_H$ ($\omega\in \Omega^\ut_{I;A}$). 

\smnoind
(c) If all $A_i$ are $C^*$-algebras and all $\Psi_i$ are irreducible, then $\alpha^\Psi$ is transitive. 
\end{rem}

\medskip

\begin{cor}\label{cor:rep-on-inf-dir-prod}
There is a unital $^*$-representation $\prod \tsi \Psi_i: \BOI^\ut A_i \to \CL( \prod \tsi H_i )$ such that 
for any $\mu\in \Omega^\un_{I;H}$, $\omega\in \Omega^\ut_{I;A}$ and $b\in \BOI^\omega A_i$, 
\begin{equation}\label{eqt:PIti<->BOI}
\big(\prod \tsi \Psi_i\big)(b)\circ \ti\Upsilon^{\mu} 
\ = \ \ti\Upsilon^{\omega\cdot \mu} \circ \big(\BOI^{\phi_1}\Psi_i\big)(b)\big|_{\bar\bigotimes_{i\in I}^\mu H_i}, 
\end{equation}
where $\ti\Upsilon^\mu$ is as in Proposition \ref{prop:cp-vN}(b). 
\end{cor}
\begin{prf}
By Proposition \ref{prop:cp-vN}(b), 
there is a bounded linear map 
$$\big(\prod \tsi \Psi_i\big)(b): \prod \otimes_\ii^{\kappa_H(\mu)} H_i \to \prod \otimes_\ii^{\omega\cdot\kappa_H(\mu)} H_i$$ 
such that Equality \eqref{eqt:PIti<->BOI} holds (see also Remark \ref{rem-act-OA-OH}(b)). 
Since we have $\sup_{\mu\in \Omega^\un_{I;H}} \big\|(\BOI^{\phi_1}\Psi_i)(b)|_{\bar\bigotimes_{i\in I}^\mu H_i}\big\| < \infty$ (because of Theorem \ref{thm:inf-ten-c-st-alg}(a)), we know from Proposition \ref{prop:cp-vN}(a) and \cite[Lemma 4.1.1]{vN} that $(\prod \tsi \Psi_i)(b)$ induces an element in $\CL(\prod \tsi H_i )$, which clearly gives a $^*$-representation. 
\end{prf}

\medskip

It is natural to ask if $\prod \tsi \Psi_i$ is injective if all $\Psi_i$ are. 
However, $\prod \tsi \Psi_i$ is never injective as can be seen in 
Example \ref{eg:ten-prod-rep-BC}(b) and the discussion following it.  

\medskip

\begin{eg}\label{eg:ten-prod-rep-BC}
For any $\ii$, let $A_i = \BC = H_i$ and $\iota_i: A_i \to \CL(H_i)$ be the canonical map. 
Suppose that $\Phi$, $\varphi$ and $\hat \Phi$ are as in Example \ref{eg:ten-Hil-sp}. 

\smnoind
(a) Let $\Lambda: \BC[\Omega^\ut_{I;\BC}] \to \CL(\ell^2(\Omega^\ut_{I;\BC}))$ be the left regular representation. 
For every $\alpha,\beta\in \PI \BT$, one has 
$$\big(\hat\Phi^*\circ \Lambda(\lambda_{[\alpha]_\sim})\circ  \hat\Phi\big)(\tsi \beta_i)
\ = \ \varphi(\alpha^{-1}) \tsi \alpha_i\beta_i
\ = \ \big(\BOI^{\phi_1} \iota_i\big)(\Phi^{-1}(\lambda_{[\alpha]_\sim}))(\tsi \beta_i).$$
Consequently, $\BOI^{\phi_1} \iota_i$ can be identified with $\Lambda$ (under $\Phi$ and $\hat \Phi$). 

\smnoind
(b) Let $\alpha\in \PI\BT$ such that $\alpha\nsim 1$ but $\alpha \approx 1$ with $\PI \alpha_i =1$. 
If $\beta\in \PI \BC$ is a $C_0$-sequence  with $\|\prod \tsi \beta_i\| = 1$, one has $\|\prod \tsi \alpha_i\beta_i\| = 1$ and 
$$\big\langle \prod \tsi \alpha_i\beta_i, \prod \tsi \beta_i \big\rangle 
\ =\ 1,$$
which imply that $\prod \tsi \alpha_i\beta_i = \prod \tsi \beta_i$. 
Therefore, $(\prod \tsi \iota_i)(\tsi \alpha_i) = \id$ but $\tsi \alpha_i \neq \tsi 1$. 
Consequently, $\prod \tsi \iota_i$ is non-injective (actually, $(\prod \tsi \iota_i)\circ \Phi^{-1}$ is non-injective as a group representation of $\Omega^\ut_{I;\BC}$). 
\end{eg}

\medskip

In general, even $\big(\prod \tsi \Psi_i\big)|_{\BOI^\ut \BC e_i}$ is non-injective. 
In fact, suppose that $\alpha$ is as in the above. 
For any $C_0$-sequence $\xi\in \PI H_i$, with $\|\prod \tsi \xi_i\| = 1$, the same argument as 
Example \ref{eg:ten-prod-rep-BC}(b) tells us that $\prod \tsi \alpha_i\xi_i = \prod \tsi \xi_i$. 
Thus, $\big(\prod \tsi \Psi_i\big)(\tsi e_i - \tsi \alpha_i e_i) = 0$. 

\medskip

On the other hand, by Theorem \ref{thm:inf-ten-c-st-alg} and Corollary \ref{cor:rep-on-inf-dir-prod}, there exist canonical $^*$-homomorphisms  
$$J^{\phi_1}: \BOI^\ut \CL(H_i) \to \CL\big(\bar\bigotimes_\ii^{\phi_1} H_i\big)
\ \ \text{and} \ \ 
J^\Pi: \BOI^\ut \CL(H_i) \to \CL\big(\prod \tsi H_i\big).$$
Notice that $J^{\phi_1}$ is injective but $J^\Pi$ is never injective. 

\medskip

\begin{cor}\label{cor:ten-rep-prod-gps}
Let $\pi_i: G_i\to U_{\CL(H_i)}$ be a unitary representation of a group $G_i$, for each $i\in I$. 

\smnoind
(a) There exist canonical unitary representations $\BOI^{\phi_1} \pi_i$ and $\prod \tsi \pi_i$ of $\PI G_i$ on $\bar\bigotimes_\ii^{\phi_1} H_i$ and $\prod \tsi H_i$ respectively. 

\smnoind
(b) If the induced $^*$-representation $\hat \pi_i: \BC[G_i]\to \CL(H_i)$ is injective for all $i\in I$, the induced $^*$-representation $\widehat{\BOI^{\phi_1} \pi_i}$ of $\BC[\PI G_i]$ is also injective. 
\end{cor}
\begin{prf}
(a) Let $\BOI^\ut \pi_i := \Theta_{\CL(H)}\circ \PI \pi_i: \PI G_i \to \BOI^\ut \CL(H_i)$. 
Then 
$$\BOI^{\phi_1} \pi_i\ :=\ J^{\phi_1}\circ \BOI^\ut \pi_i 
\quad \text{and}\quad 
\prod\tsi \pi_i\ :=\ J^{\Pi}\circ \BOI^\ut \pi_i$$ 
are the required representations. 

\smnoind
(b) By Theorem \ref{thm:inf-ten-c-st-alg}(c), $\BOI^{\phi_1} \hat\pi_i$ is injective. 
As $\widehat{\BOI^{\phi_1} \pi_i}$ is the restriction of $\BOI^{\phi_1} \hat\pi_i$ on $\BC[\PI G_i]$ (see Example \ref{eg:gp-alg-inf-prod-gp}(a)), it is also injective. 
\end{prf}

\medskip

\begin{cor}
$\prod \tsi \Psi_i$ is never irreducible, and neither do $\BOI^{\phi_1}\Psi_i$. 
\end{cor}
\begin{prf}
Let $\tau_i: \BC\to A_i$ be the canonical unital map and set $\check \Psi_i := \Psi_i\circ \tau_i$ ($i\in I$).
Suppose that $\alpha, \beta\in \PI \BT$ with $\alpha \not\approx \beta$ and $\xi\in \PI^\un H_i$. 
Then $[\alpha_i\xi_i]_\ii \not\approx [\beta_i\xi_i]_\ii$ and the two unit vectors 
$$\big(\prod \tsi\check\Psi_i\big)(\tsi \alpha_i)\big(\prod \tsi \xi_i\big)
\quad \text{and} \quad 
\big(\prod \tsi\check\Psi_i\big)(\tsi \beta_i)\big(\prod \tsi \xi_i\big)$$ 
are orthogonal. 
Consequently, $\dim\ \!(\prod \tsi \check \Psi_i)(\BC^{\otimes I}_\ut) > 1$. 
As $(\prod \tsi \Psi_i)\circ (\BOI \tau_i) = \prod \tsi \check\Psi_i$, we have $(\prod \tsi \check\Psi_i)(\BC^{\otimes I}_\ut) \subseteq Z\big((\prod \tsi \Psi_i)(\BOI^\ut A_i)\big)$ and $\prod \tsi \Psi_i$ is not irreducible.
A similar but easier argument also shows that $\BOI^{\phi_1}\Psi_i$ is not irreducible. 
\end{prf}

\medskip

For any $C^*$-algebra $A$, we denote by $S(A)$ the state space of $A$ and by $(H_\rho, \pi_\omega, \xi_\omega)$ the GNS construction of $\omega\in S(A)$. 
We would like to consider a natural injective $^*$-representation of $\BOI^\ut A_i$ defined in terms of $(H_{\omega_i}, \pi_{\omega_i})$.  

\medskip

If $\rho\in \PI  S(A_i)$ and $\check \rho$ is defined as
$$\check \rho(a) := \big\langle \big(\BOI^{\phi_0} \pi_{\rho_i}\big)(a)(\tsi \xi_{\rho_i}), (\tsi \xi_{\rho_i})\big\rangle 
\qquad \big(a\in \BOI^\ut A_i\big),$$
then the closure of $\big(\BOI^{\phi_1} \pi_{\rho_i}\big)(\BOI^\ut A_i)(\tsi \xi_{\rho_i})$ will coincide with $H_{\check \rho} := \bar\bigoplus_{\omega\in \Omega^\ut_{I;A}}\bar\bigotimes_{i\in I}^{\omega\cdot [\xi_\rho]_\sim} H_{\rho_i} \subseteq \bar\bigotimes_\ii^{\phi_1} H_{\rho_i}$.
We set 
$\pi_{\check \rho}(a) := \big(\BOI^{\phi_1} \pi_{\rho_i}\big)(a)|_{H_{\check\rho}}$. 
Notice that if all $\rho_i$ are pure states, then $H_{\check\rho} = \bar\bigotimes_\ii^{\phi_1} H_{\rho_i}$ (see Remark \ref{rem-act-OA-OH}(c)). 

\medskip

\begin{cor}\label{cor:spat-ten-prod}
Let $A_i$ be a $C^*$-algebra ($\ii$).
The $^*$-representation $\Psi_A:=\bigoplus_{\rho\in \PI S(A_i)} (H_{\check \rho}, \pi_{\check\rho})$ is injective.
Consequently, the $^*$-representation $\Phi_A:=\bigoplus_{\rho\in \PI S(A_i)} \big( {\bar\bigotimes}_{i\in I}^{\phi_1} H_{\rho_i}, \BOI^{\phi_1} \pi_{\rho_i}\big)$ is also injective. 
\end{cor}
\begin{prf}
Suppose that $(H_i, \Psi_i)$ is a universal $^*$-representation of $A_i$ ($i\in I$). 
Let $F, u^{(1)}, ..., u^{(n)}, b^{(1)},...,b^{(n)}$ as well as $a^{(1)},...,a^{(n)}$ be as in the proof of Theorem \ref{thm:inf-ten-c-st-alg}(c) with  
$\Psi_A\Big({\sum}_{k=1}^n a^{(k)}\Big) = 0$.
Again, it suffices to show that $a^{(1)} = 0$, and we may assume that $u^{(1)}_i = e_i$ ($i\in I$).
If $n = 1$, we take any $x\in \PI  \sph(H_i)$. 
If $n>1$, we take an element $x\in \PI  \sph(H_i)$ satisfying
$$x\ \nsim\ \big[\Psi_i\big(u^{(k)}_i\big)x_i\big]_{i\in I} 
\qquad (k = 2,...,n)$$
(the argument of Theorem \ref{thm:inf-ten-c-st-alg}(c) ensures its existence). 
Let us set 
$\rho_i(a):= \langle \Psi_i(a)x_i,x_i\rangle$
when $i\in I\setminus F$, 
and pick any $\rho_i\in S(A_i)$ when $i\in F$. 
For every $i\in I\setminus F$, one may regard $\big(H_{\rho_i}, \pi_{\rho_i}\big)$ as a subrepresentation of $(H_i, \Psi_i)$ such that $\xi_{\rho_i}\in H_{\rho_i}$ is identified with $x_i\in H_i$. 
Then $x$ can be considered as an element in $H_{\check\rho}$. 
Since $x \nsim \big[\pi_{\rho_i}\big(u^{(k)}_i\big)x_i\big]_{i\in I}$ for all $2\leq k \leq n$, the argument of Theorem \ref{thm:inf-ten-c-st-alg}(c) tells us that 
$$\big(\BOI^{[x]_\sim} \pi_{\rho_i}\big)(a^{(1)})\eta
\ =\ 0 
\qquad \big(\eta\in \BOI^x H_{\rho_i}\big).$$ 
Consequently, $\big(\BOF\pi_{\rho_i}\big)\big(b^{(1)}\big) = 0$ and 
$b^{(1)} =0$ (as $\rho_i$ is arbitrary when $i\in F$).
The second statement follows readily from the first one. 
\end{prf}

\medskip

Notice also that $\big( {\bar\bigotimes}_{i\in I}^{\phi_1} H_{\rho_i}, \BOI^{\phi_1} \pi_{\rho_i}\big)$ is in general not a cyclic representation, and $(H_{\check \rho}, \pi_{\check\rho})$ can be regarded as a cyclic analogue of it.

\medskip

We end this paper with the following result concerning tensor product of Hilbert algebras. 

\medskip

\begin{cor}\label{cor:hil-alg}
Let $\{A_i\}_\ii$ is a family of unital Hilbert algebras (see e.g. \cite[Definition VI.1.1]{Take}) such that $\|e_i\| =1$ ($\ii$). 
Then $A:=\BOI^\ut A_i$ is also a unital Hilbert algebra with $\|\tsi e_i\| = 1$. 
\end{cor}
\begin{prf}
Note that since $\| e_i \| =1$, one has $\|u_i\| = 1$ for any $u_i\in U_{A_i}$. 
Thus, we have $\BOI^\ut A_i \subseteq \BOI^\un A_i$, which gives an inner product $\langle \cdot, \cdot\rangle_A$ on $A$. 
Observe that $\BOI^{\omega} A_i$ is orthogonal to $\BOI^{\omega'} A_i$ (in terms of $\langle \cdot, \cdot\rangle_A$) whenever $\omega$ and $\omega'$ are distinct elements in $\Omega^\ut_{I;A}$.
Thus, in order to show the involution on $A$ being an isometry, it suffices to check that $\|x^*\| = \|x\|$ whenever $x\in \BOI^\omega A_i$ and $\omega\in \Omega^\ut_{I;A}$. 
In fact, for any $u\in \PI U_{A_i}$, $F\in \KF$ and $a\in {\bigotimes}_{i\in F} A_i$, we have  
$$\|J_F^u(a)^*\| 
\ = \ \|J^{u^*}_F(a^*)\| 
\ = \ \|a^*\| 
\ = \ \|a\| 
\ = \ \|J^u_F(a)\|,$$
because the involution on ${\bigotimes}_{i\in F} A_i$ is an isometry. 
Let $H_i$ be the completion of $A_i$ (with respect to the inner-product) and $\Psi_i: A_i \to \CL(H_i)$ be the canonical unital $^*$-representation ($\ii$). 
Since
$$\BOI^{\phi_1} \Psi_i(a)b
\ = \ ab
\qquad (a,b\in A),$$
Theorem \ref{thm:inf-ten-c-st-alg}(a) tells us that for each $x\in A$, one has $\langle xy,z\rangle_A = \langle y,x^*z\rangle_A$ ($y,z\in A$) and $\sup_{\|y\|\leq 1}\|xy\| < \infty$. 
Finally, as $A$ is unital, we see that $A$ is a Hilbert algebra (with $\|\tsi e_i\| =1$). 
\end{prf}

\medskip

Consequently, if all $A_i$ are weakly dense unital $^*$-subalgebras of finite von-Neumann algebras, then so is $\BOI^\ut A_i$.

\medskip\noindent
Chi-Keung Ng, Chern Institute of Mathematics, Nankai University, Tianjin 300071, China.

\smnoind
\emph{Email address:} ckng@nankai.edu.cn


\medskip

\begin{thebibliography}{9}

\bibitem{AN}
H. Araki and Y. Nakagami, A remark on an infinite tensor product of von Neumann algebras, Publ. Res. Inst. Math. Sci. 8 (1972), 363-374.

\bibitem{BeC}
E. B\'{e}dos and R. Conti, On infinite tensor products of projective unitary representations, Rocky Mountain J. Math. 34 (2004), 467-493. 

\bibitem{BC}
W. Bergmann and R. Conti, On infinite tensor products of Hilbert $C^*$-bimodules, \emph{Operator algebras and mathematical physics (Constanta 2001)}, Theta, Bucharest (2003), 23-34. 

\bibitem{Black}
B. Blackadar, Infinite tensor products of $C^*$-algebras, Pac. J. Math. 77 (1977), 313-334. 

\bibitem{Bru}
H. W. Bruce, Infinite Tensor Products of Commutative Subspace Lattices, Proc. Amer. Math. Soc. 103 (1988), 429-437.

\bibitem{BS}
R. C. Busby and H. A. Smith, Representations of twisted group algebras, Trans. Amer. Math. Soc. 149 (1970), 503-537. 

\bibitem{Cal}
L. Calabi, Sur les extensions des groupes topologiques, Ann. Mat. Pura Appl. 32 (1951), 295-370. 

\bibitem{EM}
S. Eilenberg and S. MacLane, Group extensions and homology, Ann. of Math. 43 (1942), 757-831.

\bibitem{EK}
D. E. Evans and Y. Kawahigashi, \emph{Quantum Symmetries on Operator Algebras}, Clarendon Press, Oxford (1998). 

\bibitem{Flor}
R. Floricel, Infinite tensor products of spatial product systems, Infin. Dimens. Anal. Quantum Probab. Relat. Top. 11 (2008), 447-465. 

\bibitem{HN}
Y. P. Huo and C. K. Ng, Some algebraic structures related to the quantum system with infinite degrees of freedom, Rep.\ in Math.\ Phy. 67 (2011), 97-107. 

\bibitem{Gill}
T. L. Gill, Infinite Tensor Products of Banach Spaces I, J. Funct. Anal. 30 (1978), 17-35.

\bibitem{GS}
T. Giordano and G. Skandalis, On infinite tensor products of factors of type ${\rm I}_2$, Erg. Theory Dynam. Sys. 5 (1985), 565-586. 

\bibitem{Gui}
A. Guichardet, \emph{Tensor products of $C^*$-algebra II: Infinite Tensor products}, Aarhus Univ. Lect. Note Ser. 13 (1969). 

\bibitem{GN}
H. Grundling and K.-L. Neeb, Infinite Tensor products of $C_0(\mathbb{R})$: Towards a Group Algebra for $\mathbb{R}^{(\BN)}$, preprint (2010), arXiv:1001.1012v1. 

\bibitem{Pow}
S. Power, Infinite tensor products of upper triangular matrix algebras, Math. Scand. 65 (1989), 291-307. 

\bibitem{RSW}
I. Raeburn, A. Sims and D. P. Williams, Twisted actions and obstructions in group cohomology, \emph{$C^*$-algebras (M\"{u}nster, 1999)}, Springer (2000), 161-181. 

\bibitem{Sto71}
E. Stormer, On Infinite Tensor Products of Von Neumann Algebras, Amer. J. Math. 93 (1971), 810-818. 

\bibitem{Take}
M. Takesaki, \emph{Theory of operator algebras II}, Encyclopaedia of Mathematical Sciences 125, Springer-Verlag, Berlin (2003). 

\bibitem{TW}
T. Thiemann and O. Winkler, Gauge Field Theory Coherent States (GCS): IV. Infinite Tensor Product and Thermodynamical Limit,  Classical Quantum Gravity 18 (2001), 4997-5053.

\bibitem{vN} 
J. von Neumann, On infinite direct products, Compositio Math. 6 (1937), 1-77.
\end{thebibliography}
\end{document}